\documentclass[smallcondensed]{svjour3}     

\usepackage{amssymb}
\usepackage{amsmath}
\usepackage{graphicx}
\usepackage{algorithm}
\usepackage{algorithmicx}
\usepackage[noend]{algpseudocode}
\usepackage{algpseudocode}
\usepackage{color}

\journalname{Journal of Scientific Computing}%

\begin{document}
\title{A scalable space-time domain decomposition approach for solving large-scale nonlinear regularized inverse ill-posed problems in 4D variational data assimilation}

\titlerunning{Decomposition approaches for solving 4D-Var DA}        

\author{Luisa D'Amore \and 
Emil Constantinescu \and 
Luisa Carracciuolo
}


\institute{Luisa D'Amore  \and
\at
            University of Naples Federico II, Naples, Italy\\
            (\email{luisa.damore)@unina.it}  \\
            \and
            Emil Costantinescu \at Mathematics and Computer Science Division, Argonne National Laboratory, Lemont, Illinois, USA\\
            (\email{emcosta@mcs.anl.org})\\
            \and
            Luisa Carracciuolo \at Istituto per i Polimeri, Compositi e Biomateriali of the CNR (IPCB-CNR), Rome, Italy\\
            (\email{luisa.carracciuolo@cnr.it})          
}

\date{Received: date / Accepted: date}

\maketitle
\begin{abstract}
We develop innovative  algorithms for solving the strong-constraint formulation of four-dimensional variational data assimilation  in large-scale applications. We present a space-time decomposition  approach that employs domain decomposition along both the spatial and temporal directions in the overlapping case and involves  partitioning of both the solution and the operators. Starting from the global functional defined on the entire domain, we obtain a type of regularized local functionals on the set of subdomains providing the  order reduction of both  the predictive and the data assimilation models. We analyze the algorithm convergence and its performance in terms of reduction of time complexity and  algorithmic scalability. The numerical experiments are carried out   on the shallow water equations on the sphere 
according to the setup available at the {\em Ocean Synthesis/Reanalysis Directory} provided by Hamburg University. 
\end{abstract}

\section{Introduction and motivation}
\noindent  
\noindent 
Assimilation of observations into models is a well-established critical practice in the meteorological community. Operational models require on the order of $10^7$ or $10^8$  model variables and  the capacity to assimilate on the order of $10^6$ observations. Various approaches have been proposed for reducing the complexity of assimilation methods to make  them more computationally affordable while retaining their original accuracy. Ensemble approaches and reduced-order models are the most significant approximations.  Other approaches take full advantage of existing partial differential equations (PDEs)-based solvers,  based on spatial domain decomposition (DD) methods, where the DD solver is suitably modified
to also solve the  adjoint associated with the forward model. A different approach is the combination of DD methods in space and data assimilation (DA), where a spatial
domain-decomposed uncertainty quantification approach performs DA at the local level by using Monte
Carlo sampling \cite{Antil2010,Amaral,Liao}.
The parallel data assimilation framework \cite{PDAF} implements parallel ensemble-based Kalman filters  coupled
with the PDE-model solver.\\ 
\noindent These  methods reduce the spatial dimensionality of the predictive model, and  the resulting reduced-order model  is then resolved in time via numerical integration, typically
with the same time integrator and time step employed for the high-fidelity model leading to high-precision time synchronization.   
In the past decades, parallel-in-time methods have been investigated for reducing the temporal dimensionality of  evolutionary problems. 
Pioneering work includes that of Nievergelt  (1964),  who proposed the first time decomposition algorithm for finding the parallel solutions of evolutionary ordinary differential equations, and that of Hackbusch (1984),  who noted  that relaxation operators in multigrid can be employed on multiple
time steps simultaneously. Since then, time parallel time integration methods have been extensively expanded. A large literature list can be found at \cite{PINT-site}, which collects information about the community, methods, and software in the field of parallel-in-time integration methods. Recent efforts include
the parallel full approximation scheme in space and time (PFASST), introduced by \cite{Emmett}. PFASST  reduces the
optimization overhead by integrating the PDE-based model  directly into the optimization process,
thus solving the  PDE, the adjoint equations, and the optimization
problem simultaneously. A nonintrusive framework for integrating existing unsteady PDE solvers into a parallel-in-time simultaneous optimization algorithm, using PFASST,  is provided in  \cite{Gunther}.  Related  parallel PDE solvers based on a Schwarz preconditioner  in
space-time are  proposed in \cite{Gander2016,Liu,Ulbriq}. 

\noindent In this study we present the design of an innovative mathematical model and the development and analysis of the related numerical algorithms,  based on the simultaneous introduction of  space-time decomposition  in the overlapping case on  the PDEs  governing the physical model and  on the DA model. The core of our approach  is that the DA model acts as coarse predictor operator solving the local PDE model, by providing the background values as initial conditions of the local PDE models. Moreover, in contrast to the other  decomposition-in-time approaches, in our approach local solvers (i.e., both the coarse and the fine solvers)  run concurrently from the beginning.  Consequently, the resulting algorithm requires only the exchange of boundary conditions between adjacent subdomains. The proposed  method  belongs to the so-called reduced-space optimization techniques, in contrast to  full-space approaches such as  the PFASST method, reducing   the
runtime of the forward and the backward integration time loops. Consequently,  we could combine the proposed  approach with  the PFASST algorithm. Indeed, PFASST could be concurrently employed as the local solver of  each reduced-space  PDE-constrained optimization subproblem, exposing even more temporal parallelism.\\

\noindent Specific contributions of this work include (1) a novel decomposition approach in space-time leading to a reduced-order model of the coupled PDE-based 4D-Var DA problem; (2) strategies for computing the ``kernels'' of the resulting regularized  nonlinear least squares computational problem; and (3) a priori  performance analysis that enables a suitable implementation of the  algorithm in advanced computing environments. Results presented here are  intended as the starting point for the software development to make decisions about computer architecture, future estimates of the problem size (e.g., the resolution of the model and the number of observations to be assimilated), and  the  performance  and  parallel  scalability  of  the  algorithms.



\noindent The article is organized as follows. Section 2 gives a brief introduction to the data assimilation framework, where we follow the  discretize-then-optimize approach. The main result is the 4D-Var functional decomposition, which is given in Section 3.  In  Section 4 we review the whole parallel algorithm; its performance analysis  is discussed  in Section 5 on the shallow water equations on the sphere. The number of state variables in the
model, the number of observations in an assimilation
cycle, and the numerical parameters as the discretization step in the time and   space domains are defined on the basis of a discretization grid using data from the Ocean Synthesis/Reanalysis Directory of Hamburg University (\cite{Dati}).  A scalability prediction of the case study based on the shallow water equations is presented in Section 6. Our conclusions are provided in Section 7.  

\section{Data assimilation framework}
We begin with a  general DA problem setup and then  consider a more convenient setup for describing the domain decomposition approach. \\

\noindent Let  $\mathcal{M}^{\Delta\times \Omega}$ denote  a forecast model described by  nonlinear Navier--Stokes equations,\footnote{Examples are the primitive equations of oceanic circulation models that are based on Boussinesq, hydrostatic momentum, mass balances, material tracer conservation, the seawater equation of state, and parameterized subgrid-scale transports \cite{MooreI,MooreII,MooreIII,Moore,NEMO}.} where $\Delta \subset \Re$ is the  time interval and $\Omega \subset \Re^N$ is the spatial domain. If  $t \in \Delta$ denotes the time variable and  $x \in \Omega$ the spatial variable, let\footnote{Although typical prognostic variables are temperature, salinity,
horizontal velocity, and sea surface displacement, here, for simplicity of notations,  we assume that $u^b(t,x) \in \Re$.} $$u^b(t,x) : \, \Delta  \times \Omega \mapsto \Re$$ is the function  representing the solution of  $\mathcal{M}^{\Delta\times \Omega}$, which we assume  belongs to the Hilbert space $\mathcal{K}(\Delta \times \Omega)$ equipped with the standard Euclidean norm. Following \cite{Daget}, we assume that  $\mathcal{M}^{\Delta \times \Omega}$ is symbolically described  as the following initial value problem:

\begin{equation}\label{cond1}
\left \{
 \begin{array}{ll}
  u^b(t,x)= \mathcal{M}^{\Delta \times \Omega}[u^b(t_0,x)],  & \forall \, (t,x) \in \Delta \times \Omega, \\
  u^b(t_0, x)=u_0^b(x),  & t_0\in \Delta \, , \,x \in \Omega \,. \\
 \end{array}
\right .
\end{equation}

\noindent The function $u^b(t,x)$ is referred to as the background state in $\Delta \times \Omega$. The function $u^b_0(x)$ is the initial condition of $\mathcal{M}^{\Delta \times \Omega}$, and this is the value of the background state in $t_0\times \Omega$. Let

\begin{equation}\label{cond2}
  v(\tau,y)=\mathcal{H}(u(t,x)), \quad (t,x)  \in \Delta \times  \Omega, \quad (\tau,y)  \in \Delta' \times  \Omega' ,
  \end{equation}
\noindent where $\Delta'\subset \Delta$ is the observation time interval and $\Omega' \subset \Re^{nobs}$, with $\Omega'\subset \Omega$, is the observation spatial domain.  $$\mathcal{H}: \mathcal{K}(\Delta \times \Omega) \mapsto \mathcal{K}(\Delta'\times \Omega')$$ denotes the observation mapping, where $ \mathcal{H}$ is a nonlinear operator that includes transformations and grid interpolations.\\
According to the practical applications  of model-based assimilation of observations, we  use the following definition of a data assimilation  problem associated with $\mathcal{M}^{\Delta \times \Omega}$.

\begin{definition}[DA problem setup]
We consider the following setup.\footnote{Throughout the paper, for simplicity, we use the notation $j=1,K$ to indicate  $j=1, \ldots, K$.}
\begin{itemize}
  \item Let $\{t_k\}_{k=0,M-1}$,  where $t_k=t_0+k \Delta t$,  be  a discretization of $\Delta$,   such that
$\Delta_M:=[t_0,t_{M-1}] \subseteq \Delta \,.$
  \item Let $D_{K}(\Omega):=\{x_j\}_{j=1,K}\in \Re^{K}$,  be  a discretization of  $\Omega$  such that $D_{K}(\Omega) \subseteq \Omega\,.$
  \item  $\Delta_M\times \Omega_{K}=\{\mathbf{z}_{ji}:=(t_j,x_i)\}_{i=1,K;j=1,M}$.
  \item Let $\mathbf{u}_0^{b} :=\{u_0^j\}_{j=1,K}^{b} \equiv \{ u(t_0,x_j)^{b} \}_{j=1,K}\in \Re^{K}$ be the discretization of initial value in (\ref{cond1}).
  \item Let $\mathbf{u}^b= \{\mathbf{u}_k^b\}_{k=0,M-1}$ where  $\mathbf{u}_k^b:= \{ u^b(t_k,x_j) \}_{j=1,K}\in \Re^K$ be the  numerical solution   of  (\ref{cond1}) at $t_k$.
      \item Let $nobs<<K$.
      \item  Let $\Delta'_M=[\tau_0, \tau_{M-1}]\subseteq    \Delta_M$.
\item Let $D'_{nobs}(\Omega'):=\{x_j\}_{j=1,nobs}\in \Re^{nobs}$  be  a discretization of  $\Omega'$  such that $D_{nobs}(\Omega') \subseteq \Omega'$\, .
  \item Let $\mathbf{v}= \{\mathbf{v}_k\}_{k=0,M-1}$  where $\mathbf{v}_k:=\{v(\tau_k,x_j)\}_{j=1,nobs}\in \Re^{nobs}$ be  the  values of the observations on $x_j$ at $\tau_k$.
      \item Let $\{\mathbf{H}^{(k)}\}_{k=0,M-1}$, the tangent linear model (TLM)  of $\mathcal{H}(u(t_k,x))$ at time $t_k$.
      \item Let $\mathbf{M}^{\Delta_M \times \Omega_{K}}$ be a discretization of  $\mathcal{M}^{\Delta\times \Omega}$.
          \item Let $\mathbf{M}^T$ is the adjoint  model (ADM)\footnote{Let $\mathbf{A}:\mathbf{x}\to \mathbf{y}=\mathbf{A}\mathbf{x} $ be a linear operator on $\Re^{N}$ equipped with the standard Euclidean norm.
The operator $\mathbf{A}^T:\mathbf{y}\to
\mathbf{x}=\mathbf{A}^T\mathbf{y} $, such that 
\begin{equation}\label{Eq_propAdj}
<\mathbf{y}, \mathbf{A}\mathbf{x}>=<\mathbf{A}^T\mathbf{y},\mathbf{x}>, \quad \forall \mathbf{x},\forall \mathbf{y},
\end{equation}
where $<\cdot,\cdot>$ denotes the scalar product in $\Re^N$, is the adjoint  of $\mathbf{A}$.}
of  $\mathbf{M}^{0,M-1}$
\cite{GIERING}\footnote{If $\mathbf{M}^{i-1,i}$ is the TLM of $\mathcal{M}^{\Delta\times \Omega}$, in $[t_{i-1},t_i]\times \Omega_{K}$ ,  then it holds that
\begin{equation}\label{adjProd}
(\mathbf{M}^{0,M-1})^T=(\mathbf{M}^{0,1}\cdot\mathbf{M}^{1,2}\cdots\mathbf{M}^{M-2,M-1})^T=(\mathbf{M}^{M-2,M-1})^T\cdots (\mathbf{M}^{1,2})^T  (\mathbf{M}^{0,1})^T.
\end{equation}}.

\end{itemize}
\begin{flushright}
$\spadesuit$
\end{flushright}
\end{definition}
The aim of  DA is to produce the optimal combination of the
background and observations throughout the assimilation window $\Delta'_M$, in other words,  to find an optimal  tradeoff between the  estimate of the system state $\mathbf{u}^b$ and  $\mathbf{v}$.  The best estimate that optimally
fuses all this  information is called the analysis, and it is denoted as $\mathbf{u}^{DA}$. It is then used as an initial condition for the next forecast.\\

\begin{definition}[The 4D-Var DA problem: a regularized nonlinear least squares problem (RNL-LS)]
Given the DA problem setup, the 4D-Var DA problem consists of computing the vector $\mathbf{u}^{DA}\in \Re^{K}$ such that
\begin{equation}\label{min3}
\mathbf{u}^{DA} = \arg\min_{\mathbf{u} \in \Re^{K}} J(\mathbf{u})
\end{equation}
with
\begin{equation}\label{min3sec}
 J(\mathbf{u})= \| \mathbf{u}- \mathbf{u}^{b}_0\|_{{\bf B}^{-1}}^{2}    + \lambda \sum_{k=0}^{M-1} \| \mathbf{H}^{(k)}(\mathbf{M}^{\Delta_M \times \Omega_{K}} ( \mathbf{u}))-\mathbf{v}_k \|_{{\bf R}_k^{-1}}^{2},
\end{equation}

\noindent where  $\lambda>0$ is the regularization parameter;  ${\bf B}$ and ${\bf R}_k$ ($\forall k= 0,\dots,M-1$) are the covariance matrices of the errors on the background and the observations, respectively; and  $\| \cdot \|_{{\bf B}^{-1}}$ and $\| \cdot \|_{{\bf R}_k^{-1}}$ denote the weighted Euclidean norm, respectively.
\begin{flushright}
$\spadesuit$
\end{flushright}
\end{definition}

\noindent The first term in (\ref{min3sec}) quantifies the departure of the solution $\mathbf{u}^{DA}$ from
the background state $\mathbf{u}^{b}$. The second term measures the
 mismatch between the new  trajectory and observations
$\mathbf{v}_k$ for each time $t_k$  in the assimilation window. The weighting matrices ${\bf B}$
and ${\bf R}_k$ need to be predefined, and their quality influences the accuracy of the resulting analysis \cite{JCP2017}.\\

\noindent This nonlinear least-squares problem is typically considered large scale with $K$ larger than $10^6$. We next provide a mathematical formulation of a domain decomposition approach that starts from the decomposition of the whole domain $\Delta \times \Omega$ (i.e., in both space and time); it uses a partitioning of the  solution and a modified  functional describing the RNL-LS problem on the subdomain of the decomposition. Solution continuity
equations across interval boundaries are
added as constraints of the assimilation functional. We first introduce the  domain decomposition of $\Delta \times \Omega$ and then define the restriction and extension operators on functions given on $\Delta \times \Omega$. These definitions are then generalized to $\Delta_M \times \Omega_{K}$.

\section{The space-time decomposition}

In this section we give a precise mathematical setting for space and operator decomposition. 
In particular, we introduce the
functional and domain decomposition. Then, by using restriction and extension operators, we associate with the domain decomposition a functional decomposition. To this end, we  prove the following result: the minimum of the global functional, defined on the
entire domain, can be  obtained by collecting the minimum of each local functional.\\
\subsection{The space-time decomposition of the continuous 4D-Var DA model}
\noindent For simplicity we assume that the spatial  and temporal domains of the observations are the same as the background state, namely, $\Delta'=\Delta$ and $\Omega'=\Omega$; furthermore, we assume that $t_k=\tau_k$.   
\begin{definition}[Domain decomposition]\label{domdecom}
Let $P\in \mathbf{N}$ and $Q\in\mathbf{N}$ be fixed. The set of bounded Lipschitz domains ${\Omega_i}$,  overlapping subdomains of $\Omega$,
\begin{equation} \label{DDOmega}
DD(\Omega)=\left\{ \Omega_i \right\}_{i=1,P},
\end{equation}
 is called a decomposition of
$\Omega$ if
\begin{equation}\label{decom}
\bigcup_{i=1}^{P} \Omega_i =\Omega
\end{equation}
\noindent  with $$ \Omega_{jh}:=\Omega_j \cap \Omega_h\neq \emptyset$$
when two subdomains are adjacent. Similarly,  the set of overlapping subdomains of $\Delta$,
\begin{equation} \label{DDOmegaT}
DD(\Delta)=\left\{ \Delta_j \right\}_{j=1,Q},
\end{equation}
 is a decomposition of
$\Delta$ if
\begin{equation}\label{decomT}
\bigcup_{j=1}^{Q} \Delta_j =\Delta
\end{equation}
\noindent  with $$\Delta_{ik}:=\Delta_i\cap \Delta_k \neq \emptyset$$ when the two subdomains are adjacent.
\noindent We denote the  domain decomposition of $\Delta \times \Omega$ by  $DD(\Delta \times \Omega)$ with the set of $P\times Q$ overlapping subdomains of $\Delta \times \Omega$:
\begin{equation}\label{decomall}
DD(\Delta\times \Omega)= \left\{ \Delta_j \times \Omega_i \right\}_{j=1,Q;\ i=1,P}\,.
\end{equation}
\begin{flushright}
$\spadesuit$
\end{flushright}
\end{definition}

\noindent From  (\ref{decomall}) it follows that 
$$\Delta \times \Omega = \cup \Delta_j \times \cup \Omega_i = \cup(\Delta_j \times \Omega_i)\quad .$$

\noindent Next we define the restriction operator on functions in  $\mathcal{K}(\Delta\times \Omega)$ associated with the decomposition  (\ref{decomall}).

\begin{definition}[Restriction  of a function]\label{defrestrfun}
Let
$$RO_{ji}: f \in\mathcal{K}(\Delta\times \Omega) \mapsto RO_{ji}[f]\in \mathcal{K}(\Delta_j \times \Omega_i)$$ be the restriction operator (RO) of $f$ in $DD(\Delta \times \Omega)$ as in (\ref{decomall})  such that
$$ RO_{ji}[f(t,x)]\equiv
\left \{
\begin{array}{ll}
  f(t,x), & \quad \forall \,\,(t,x)  \in \Delta_j \times \Omega_i \\
  \frac{1}{2}f(t,x), & \forall \,(t,x) \,s.t.\,x  \in \Omega_i, \quad \exists \, \bar{k}\neq j: t\in \Delta_j \cap \Delta_{\bar{k}},\\
   \frac{1}{2}f(t,x), & \forall \,(t,x)\,t  \in \Delta_j, \quad \exists \, \bar{h}\neq i: x\in \Omega_i \cap \Omega_{\bar{h}},\\
   \frac{1}{4}f(t,x), &    \exists \,(\bar{h},\bar{k}) \neq (j,i): (t,x) \in (\Delta_j \cap \Delta_{\bar{h}}) \times (\Omega_i \cap \Omega_{\bar{k}}),\\
\end{array}
\right .
$$
We define
$$ f^{RO}_{ji}(t,x) \equiv RO_{ji}[f(t,x)] \quad . $$
\begin{flushright}
$\spadesuit$
\end{flushright}
\end{definition}
For simplicity, if $i \equiv j$, we denote $RO_{ii}=RO_{i}$. \\

\noindent In line with this, given a set of $Q\times P$ functions $g_{ji}$, $j=1,\dots,Q$, $i=1,\dots,P$ each  in $\mathcal{K}(\Delta_j\times \Omega_i)$, we define the  extension operator of $g_{ji}$.
\begin{definition}[Extension  of a function]
Let
 $$EO: g_{ji} \in \mathcal{K}(\Delta_j\times \Omega_i) \mapsto EO[g_{ji}] \in \mathcal{K}(\Delta \times \Omega)$$ be the extension operator (EO) of $g_{ji}$  in $DD(\Delta \times \Omega)$ as in (\ref{decomall})  such that
$$ EO[(g_{ji}(t,x)]= \left \{
\begin{array}{ll}
  g_{ji}(t,x) & \forall \,\, (t,x)\in \Delta_j \times \Omega_i ,\\
  0 & \textnormal{elsewhere}
\end{array}
\right .
$$
We define
$$ g^{EO}_{ji}(t,x) \equiv EO[g_{ji}(t,x)]\, .$$
\begin{flushright}
$\spadesuit$
\end{flushright}
\end{definition}

\noindent For any function $u\in \mathcal{K}(\Delta\times \Omega)$,  associated with the decomposition (\ref{decom}), it holds that
\begin{equation}\label{rest}
    u(t,x)= \sum_{i=1,P; j=1,Q}EO \left [u^{RO}_{ji}(t,x) \right ].
\end{equation}
\vspace*{1cm}

\noindent
Given $P\times Q$ functions $u_{ji}(t,x) \in \mathcal{K}(\Delta_i\times \Omega_j)$, the summation
\begin{equation}\label{exte}
    \sum_{i=1,P;j=1,Q}u^{EO}_{ji}(t,x)
\end{equation}
defines a function $u \in \mathcal{K}(\Delta \times \Omega)$ such that
\begin{equation}\label{Eq_ROprop}
RO_{ji}[u(t,x)]=RO_{ji} \left [\sum_{i=1,P;j=1,Q}u^{EO}_{ji}(t,x) \right ]= u_{ji}(t,x).
\end{equation}

\noindent \noindent The main outcome of this framework is  the definition of the operator $RO_{ji}$
for the 4DVar functional  defined in (\ref{min3sec}). This definition  originates  from  the definition of the restriction operator of $\mathcal{M}^{\Delta \times \Omega}$ in (\ref{cond1}), given as follows.

\begin{definition}[Restriction of  $\mathcal{M}^{\Delta \times \Omega}$]
If $\mathcal{M}^{\Delta \times \Omega}$ is  defined in (\ref{cond1}), we introduce the model $\mathcal{M}^{\Delta_j \times \Omega_i}$ to be  the  \emph{restriction}  of $\mathcal{M}^{\Delta \times \Omega}$:
 $$RO_{ji}: \mathcal{M}^{\Delta \times \Omega}(t,x)[u(t_0,x)]  \mapsto RO_{ji}[\mathcal{M}^{\Delta \times \Omega}[u(t_0,x)]]$$  defined in $\Delta_j \times \Omega_i$ such that 
\begin{equation} \label{restricted}
\left \{
\begin{array}{ll} 
 u^b(t,x)= \mathcal{M}^{\Delta_j \times \Omega_i}[u^b(t_j,x)] & \forall \,\, (t,x)\in \Delta_j \times \Omega_i \\
 u^b(t_j,x)= u^b_j(x) & t_j \in \Delta_j\, , \quad x \in \Omega_i
\end{array}
\right. .
\end{equation}

\begin{flushright}
$\spadesuit$
\end{flushright}
\end{definition}
We note that  the initial condition  $u^b_j(x)$ is  the value in $t_j$ of the solution of  $\mathcal{M}^{\Delta \times \Omega}[u(t_0,x)]$ defined in (\ref{cond1}).\\

\subsection{Space-time decomposition of the discrete model}
\noindent  Assume that  $\Delta_M\times \Omega_{K}$ can be decomposed into a sequence of $P\times Q$ overlapping subdomains $\Delta_j \times \Omega_i$ such that
$$\Delta_M \times \Omega_{K}=\bigcup_{i=1,P; \ j=1,Q} \Delta_j \times \Omega_i,$$
where
$\Omega_i \subset\Re^{r_i}$ with  $r_i \leq K$ and $\Delta_j \subset \Re^{s_j}$ with $s_j\leq M$. Moreover, assume that 
$$ \Delta_j:= [t_j, t_{j+s_j}] \, .$$ 

\begin{definition}[Restriction of the covariance matrix]
Let $\mathbf{C}(\mathbf{w})\in \Re^{K\times K}$ be the  covariance matrix of a random vector $\mathbf{w}=(w_1, w_2, \ldots,w_{K}) \in \Re^{K}$. That is,  the coefficient  $c_{i,j}$ of $\mathbf{C}$ is $c_{i,j}=\sigma_{ij} \equiv  Cov(w_i,w_j)$. With $s<K$, we define the restriction operator $RO_{st}$ onto $\mathbf{C}(\mathbf{w})$  as follows:
$$RO_{st}: \mathbf{C}(\mathbf{w})\in \Re^{K\times K} \mapsto RO_{st}[\mathbf{C}(\mathbf{w})]:=\mathbf{C}(\mathbf{w^{RO_{st}}})\in \Re^{s\times s}\,,$$
in other words,  the covariance matrix defined on  $\mathbf{w^{RO_{st}}}$.
\begin{flushright}
$\spadesuit$
\end{flushright}
\end{definition}

\noindent Hereafter, we refer to $\mathbf{C}(\mathbf{w^{RO_s}})$ using the notation $\mathbf{C_{st}}$.\\

\begin{definition}[Restriction of the operator $\mathbf{H}^{(k)}$]
We define the \emph{restriction  operator} $RO_{ji}$ of $\mathbf{H}^{(k)}$ in $DD(\Delta \times \Omega)$ as in (\ref{decomall})  as the TLM at time $t_k$ of the restriction of $\mathcal{H}$ on $\Delta_j \times \Omega_i$.
\begin{flushright}
$\spadesuit$
\end{flushright}
\end{definition}

\begin{definition}[Restriction of  $\mathbf{M}^{\Delta_M \times \Omega_{K}} $]
We let $\mathbf{M}^{\Delta_j \times \Omega_{i}}$ be the  restriction operator  $RO_{ji}$ of $\mathbf{M}^{\Delta_M \times \Omega_{K}} $ in $\Delta_j \times \Omega_i$, where
 $$RO_{ji}: \mathbf{M}^{\Delta_M \times \Omega_{K}} (\mathbf{u}_0^b)  \mapsto\mathbf{M}^{\Delta_j \times \Omega_{i}} (\mathbf{u}_0^b)=\mathbf{u}^{b}_{ji}$$  defined in $\Delta_j \times \Omega_i$.

\begin{flushright}
$\spadesuit$
\end{flushright}
\end{definition}

\begin{definition}[Restriction of the operator $\mathbf{M}^{0,M-1}$]
We define  $\mathbf{M}^{j,j+1}_i $  to be the restriction  operator $RO_{ji}$ of $\mathbf{M}^{0,M-1}$ in $DD(\Delta \times \Omega)$,  as in (\ref{decomall}). It is  the TLM of the restriction of $\mathbf{M}^{\Delta_M \times \Omega_{K}}$ on $\Delta_j \times \Omega_i$.
\begin{flushright}
$\spadesuit$
\end{flushright}
\end{definition}

\noindent With these definitions, we are now able to construct the restriction of the entire cost functional.

\begin{definition}[Restriction of 4D-Var DA]\label{Functional-restriction}
Let
$$RO_{ji}[J]: \mathbf{u}_{{ji}}\mapsto RO_{ji}[J](\mathbf{u}_{{ji}})$$
denote the restriction operator  of the 4D-Var DA functional defined in (\ref{min3sec}). It is defined as
\begin{equation}\label{funrestr}
\begin{array}{ll}
RO_{ji}[J](\mathbf{u}_{{ji}})=& \| \underbrace{RO_{ji}(\mathbf{u})}_{\mathbf{u}_{ji}}- \underbrace{RO_{ji}[\mathbf{M}^{\Delta_M \times \Omega_{K}} (\mathbf{u}^{b}_0)]}_{\mathbf{u}^{b}_{ji}}\|_{({\bf B}^{-1})_{ji}}\\ & + \lambda \sum_{k:t_k \in \Delta_j} \| \underbrace{RO_{ji}[\mathbf{H}^{(k)}]RO_{ji}[\mathbf{M}^{\Delta_M \times \Omega_{K}}  (\mathbf{u})]}_{(\mathbf{H}^{(k)})_{ji}RO_{ji}[(\mathbf{M}^{\Delta_M \times \Omega_{K}} ) (\mathbf{u}_{ji})]}-\underbrace{RO_{ji}[\mathbf{v_k}]}_{\mathbf{v}_{ji}} \|_{({\bf R}_k^{-1})_{ji}}^{2} \, .\\
\end{array}
\end{equation}

\begin{flushright}
$\spadesuit$
\end{flushright}
\end{definition}

\noindent  The local 4D-Var DA functional $ J_{ji}(\mathbf{u}_{{ji}})$ in (\ref{funrestr}) becomes
\begin{subequations}
\label{eq:cost:parallel}
\begin{align}
J_{ji}(\mathbf{u}_{{ji}})= &\underbrace{\|{\mathbf{u}_{ji}}-{\mathbf{u}^b_{ji}}\|_{(\mathbf{B}^{-1})_{ji}}}_{local \,\,state\,\,trajectory}
+\\
& \lambda \sum _{k: t_k \in \Delta_j}\underbrace{\|(\mathbf{H}^{(k)})_{ji}[\mathbf{M}^{k,k+1}_i(\mathbf{u}_{ji})]-\mathbf{v}_{ji}\|_{(\mathbf{R}_k^{-1})_{ji}}}_{local \,\,observations}.
\end{align}
\end{subequations}

\noindent In other words, the approach we are following is  first to decompose the 4D-Var functional  $J$ and then to locally linearize and solve each local functional $J_{ji}.$\\

\noindent For simplicity of notations we let 
$$ RO_{ji}[J]  \equiv J_{ \Delta_j \times \Omega_i}.$$

\noindent We note that in \eqref{funrestr}  $RO_{ji}[J](\mathbf{u}_{{ji}})$ the first term quantifies the departure of the state $\mathbf{u}_{{ji}}$ from the background state $\mathbf{u}^b_{{ji}}$ at time $t_j$ and space $x_i$. The second term measures the mismatch between the state $\mathbf{u}_{{ji}}$ and the observation $\mathbf{v}_{{ji}}$.

\begin{definition}[Extension of 4D-Var DA]
Given  $DD(\Delta \times \Omega)$ as in (\ref{decomall}), let
$$ EO[J]: J_{\Delta_j \times \Omega_i} \mapsto J^{EO}_{\Delta_j \times \Omega_i}\,\,,$$
be the  extension operator  of the 4D-Var functional defined in (\ref{min3sec}),  where
\begin{equation}\label{EOJ}
EO[J](J_{\Delta_j \times \Omega_i})= \left \{
\begin{array}{cc}
  J_{\Delta_j \times \Omega_i} & (t,x) \in \Delta_j \times \Omega_i \\
  0 & \textnormal{elsewhere}
\end{array}
\right. .
\end{equation}
\begin{flushright}
$\spadesuit$
\end{flushright}
\end{definition}

\noindent From (\ref{restJ}), it follows that the decomposition of $J$ satisfies

\begin{equation}\label{restJ}
    J\equiv \sum_{i=1,P;j=1,Q}J^{EO}_{ \Delta_j\times \Omega_i}\,\,.\\
\end{equation}

\noindent \noindent The implication in (\ref{restJ}) is that the 4D-Var problem can be defined as a set of local 4D-Var problems  as detailed in the following section.

\subsection{Local 4D-Var DA problem: the local RNL-LS problem}
Starting from the local 4D-Var functional in (\ref{eq:cost:parallel}), which is 
 obtained by applying the restriction operator to the 4D-Var functional defined  in (\ref{min3sec}), we add  a \emph{local} constraint  to the restriction. This is a type of regularization of the local 4D-Var  functional introduced  in order to enforce the continuity of each  solution of the local  problem onto the overlap region between  adjacent subdomains. The local constraint consists of the overlapping operator 
 $\mathcal{O}_{(jh)(ik)}$ defined as 
\begin{equation}\label{overl}
\mathcal{O}_{(jh)(ik)}:= \mathcal{O}_{jh}\circ\mathcal{O}_{ik},
 \end{equation}
 where the symbol $\circ$ denotes the operators composition.  Each operator in (\ref{overl}) tackles the overlapping of the solution in  the spatial dimension and  in the temporal dimension, respectively.
More precisely, for $j = 1,\dots , Q; \, i=1 ,\ldots , P$, the operator  $\mathcal{O}_{(jh)(ik)}$ represents the overlap of
the temporal subdomains $j$ and $h$ and spatial subdomains $i$ and $k$, where $h$ and $k$ are given as in Definition \ref{defrestrfun} and


\begin{equation}\label{ospatial}
  \mathcal{O}_{ik}:\mathbf{u}_{ji} \in \Delta_j \times \Omega_i  \mapsto\mathbf{u}_{(j)(ik)}\in   \Delta_j \times (\Omega_i\cap \Omega_k)
\end{equation}

\noindent and
\begin{equation}\label{otime}
\mathcal{O}_{jh}: \mathbf{u}_{(j)(ik)}\mapsto  \mathbf{u}_{(jh)(ik)}\in (\Delta_j \cap \Delta_h) \times (\Omega_i \cap \Omega_k).
\end{equation}

\begin{remark}
We observe that  in the overlapping domain $\Delta_{jh} \times \Omega_{ik}$ we get two vectors,  $\mathbf{u}_{(jh)(ik)}$, which is obtained as the restriction of $\mathbf{u}_{(ji)}= \arg\min J_{ji}(\mathbf{u}_{ji})$ to that region, and  $\mathbf{u}_{(hj)(ki)}$, which is  the restriction  of $\mathbf{u}_{(hk)}= \arg\min J_{hk}(\mathbf{u}_{hk})$ to the same region.  The order of the indexes plays a significant role from the computing perspectives.    
\end{remark}
\noindent There  are three basic cases that  we may consider in  (\ref{overl}): 
\begin{enumerate}
  \item Decomposition in space, namely, $Q=1$ and $P>1$.
Here we get  $j=Q=1$ (i.e., the  time interval is not decomposed) and $P>1$ (i.e., the spatial domain $\Omega$ is decomposed according to the domain decomposition in (\ref{decomall})).  The overlapping operator is defined as in (\ref{ospatial}). In particular, we assume that
$$\mathcal{O}_{ik}(\mathbf{u}_{ji}):=\| \underbrace{RO_{ji}(\mathbf{u}_{jk})}_{\mathbf{u}_{j(ki)}} - \underbrace{RO_{jk}(\mathbf{u}_{ji})}_{\mathbf{u}_{(j)(ik)}}\|_{(\mathbf{B}^{-1})_{ik}}.$$
  \item Decomposition in time, namely, $Q>1$ and $P=1$.
We get  $i=P=1$ (i.e., the spatial domain is not decomposed) and $Q>1$ (i.e., the time interval  is decomposed according to the domain decomposition in (\ref{decomall})). The overlapping operator is defined as in (\ref{otime}). In particular, we assume that
$$\mathcal{O}_{jh}(\mathbf{u}_{ji}):=\| \underbrace{RO_{ji}(\mathbf{u}_{hi})}_{\mathbf{u}_{(hj)i}} - \underbrace{RO_{hi}(\mathbf{u}_{ji})}_{_{\mathbf{u}_{(jh)i}}}\|_{(\mathbf{B}^{-1})_{jh}}.$$
  \item Decomposition in space-time, namely,  $Q>1$ and $P>1$. We assume that $Q>1$ and $P>1$  (i.e., both the time interval and the spatial domain are decomposed according to the domain decomposition in (\ref{decomall})).  The overlapping operator is defined as in (\ref{overl}). In particular, we assume that
$$\mathcal{O}_{(jh)(ik)}(\mathbf{u}_{ji}):=\|\mathbf{u}_{(hj)(ki)} - \underbrace{RO_{hi}(RO_{jk}(\mathbf{u}_{ji}))}_{_{\mathbf{u}_{(jh)(ik)}}}\|_{(\mathbf{B}^{-1})_{(jh)(ik)}}. $$
\end{enumerate}
We now give the new definition of the local 4D-Var DA functional.
\begin{definition}[Local 4D-Var DA]
\noindent Given $DD(\Delta \times \Omega)$ as in (\ref{decomall}),  let
\begin{eqnarray}\label{solloc1}
J_{ji}(\mathbf{u}_{{ji}}) & = & RO_{ji}[J](\mathbf{u}_{{ji}})+\mu_{ji} \ {O}_{(jh)(ik)}(\mathbf{u}_{{ji}}), 
\end{eqnarray}

\noindent where $RO_{ji}[J](\mathbf{u}_{{ji}})$ is given in (\ref{funrestr})
$\mathcal{O}_{(jh)(ik)}$,  suitably defined  on  $\Delta_{jh}\times \Omega_{ik}$, be the local 4D-Var functional.
The parameter $\mu_{ji}$ is a regularization parameter. Also let
\begin{equation}\label{solloc}
{\mathbf{u}}_{ji}^{DA} = \arg\min_{\mathbf{u}_{{ji}}}{J}_{ji}(\mathbf{u}_{{ji}})
\end{equation}
be the global minimum of ${J}_{ji}$ in $\Delta_j \times \Omega_i$.
\begin{flushright}
$\spadesuit$
\end{flushright}
\end{definition}

\noindent  More precisely, the local 4D-Var DA functional $ J_{ji}(\mathbf{u}_{{ji}})$ in (\ref{solloc1}) becomes
\begin{subequations}
\label{eq:cost:parallel1}
\begin{align}
J_{ji}(\mathbf{u}_{{ji}})= &\underbrace{\|{\mathbf{u}_{ji}}-{\mathbf{u}^b_{ji}}\|_{(\mathbf{B}^{-1})_{ji}}}_{local \,\,state\,\,trajectory}
+\\
& \lambda \sum _{k: t_k \in \Delta_j}\underbrace{\|(\mathbf{H}^{(k)})_{ji}[\mathbf{M}^{k,k+1}_i(\mathbf{u}_{ji})]-\mathbf{v}_{ji}\|_{(\mathbf{R}_k^{-1})_{ji}}}_{local \,\,observations}+\\
&\mu\underbrace{\|\mathbf{u}_{(hj)(ki)}-\mathbf{u}_{(jh)(ik)}\|_{(\mathbf{B}^{-1})_{(jk)(ih)}}}_{overlap} ,
\end{align}
\end{subequations}
where the three terms  contributing to the definition of the  local DA functional clearly come out. We note that in  (\ref{eq:cost:parallel}) the  operator $\mathbf{M}^{k,k+1}_i$, which is defined in (\ref{adjProd}), replaces $\mathcal{M}^{\Delta_j \times \Omega_i}$. \\

\noindent Next we show that the absolute minimum of operator $J$  is found among the absolute minima of \emph{local} functionals.

\subsection{Local 4D-Var DA  minimization}
\noindent  Let
\begin{equation}\label{zero0fun}
    \mathbf{\widetilde{u}_{ji}}:=({\mathbf{u}}_{ji}^{DA})^{EO}\in \Re^{M \times K} , \quad \forall \, j=1,Q;i=1,P, \
\end{equation}

\noindent where ${\mathbf{u}}_{ji}^{DA}$ is defined in (\ref{solloc}), be (the extension of) the minimum of the (global) minima of the \emph{local} functionals  ${J}_{ji}$ as in (\ref{solloc}). Let
\begin{equation}\label{zero0}
     \mathbf{\widetilde{u}^{DA}}:= {\arg\min}_{j=1,Q;i=1,P} \left\{J \left (\mathbf{\widetilde{u}_{ji}}\right ) \right \}
\end{equation}
be its minimum.

\begin{theorem} \label{theorem215}
If  $J$ is convex and $ DD( \Delta\times \Omega)$  is a decomposition of $ \Delta\times \Omega$ as defined in (\ref{decomall}), then

 \begin{equation}\label{thesis}
  J(\mathbf{{u}^{DA}}) \leq  J( \mathbf{\widetilde{u}^{DA}}),
 \end{equation}

\noindent with $ \mathbf{u}^{DA}$ defined in (\ref{min3}). 
\end{theorem}

\noindent \textbf{Proof:} Let  ${\mathbf{u}}_{ji}^{DA}$  be  defined in (\ref{solloc}); it is
\begin{equation}\label{uno}
    \nabla  {J}_{ji}[ {\mathbf{u}}_{ji}^{DA}]=\underline{0}\in \Re^{NP}, \quad \forall (j,i) : \Delta_j \times \Omega_i\in DD(\Delta \times \Omega).
\end{equation}

\noindent From (\ref{uno})  it follows that 

\begin{equation}\label{nabla-sum-EO}
   \nabla EO\left[ {J}_{ji}\left ( {\mathbf{u}}_{ji}^{DA}\right ) \right]=\underline{0},
\end{equation}

\noindent which gives from  (\ref{restJ})

\begin{equation}\label{nabla-sum-EO2}
   \nabla {J}\left [({\mathbf{u}}_{ji}^{DA})^{EO}\right ]=\underline{0}.
\end{equation}

\noindent Then $({\mathbf{u}}_{ji}^{DA})^{EO}$ is a stationary point for $J$ in $\Re^{M \times K}$. Since $ \mathbf{u}^{DA}$ in (\ref{min3}) is the global minimum of $J$ in $\Re^K$, it follows that

\begin{equation}\label{diseq-minim}
 J(\mathbf{{u}}^{DA}) \leq  J\left( (\mathbf{u}_{ji}^{DA})^{EO}\right), \quad \forall \, j=1,Q;i=1,P.
\end{equation}

\noindent Then, from  (\ref{zero0}) it follows that

\begin{equation}\label{diseq-minim2}
J( \mathbf{u}^{DA} )\leq J\left( \mathbf{\widetilde{u}^{DA}}\right ) \quad .
\end{equation}

\noindent  Now  we prove that  if $J$ is convex, then $$J( \mathbf{u}^{DA} )= J(\mathbf{\widetilde{u}^{DA}})$$ by contradiction. Assume that

\begin{equation}\label{PerAbsurd}
J(\mathbf{u}^{DA} )< J( \mathbf{\widetilde{u}^{DA}} ). 
\end{equation}

\noindent In particular, $$J(\mathbf{u}^{DA}) <J( RO_{ji}( \widetilde{\mathbf{u}}^{DA})) \quad .$$

\noindent This means  that

\begin{equation}\label{eq_ROabsurd}
 RO_{ji} \left[ J(\mathbf{u}^{DA} )\right ] <  RO_{ji}\left[ J(\mathbf{\widetilde{u}^{DA} })\right].
\end{equation}

\noindent From  (\ref{eq_ROabsurd}) and  (\ref{zero0}), it is

$$
  RO_{ji} \left[J(\mathbf{u}^{DA} )\right] <  RO_{ji}\left[ min_{ji} ( J\left(  \mathbf{u}_{ji}^{DA})^{EO} \right )  \right] .
$$

\noindent Then, from  (\ref{Eq_ROprop}):

\begin{equation}\label{Eq_theAbsurd}
J_{ji} \left(  RO_{ji} [\mathbf{u}^{DA}]^{EO} \right) < J_{ji}\left ( RO_{ji}\left[ \mathbf{u}_{ji}^{DA}\right ]^{EO}\right )  = J_{ji}({\mathbf{u}}_{ji}^{DA} )\quad \quad .
\end{equation}

\noindent Equation (\ref{Eq_theAbsurd}) is a contradiction because the value of ${\mathbf{u}}_{ji}^{DA}$ is the global minimum for $J_{ji}$, and therefore the (\ref{thesis}) is proved.

\begin{flushright}
$\clubsuit$
\end{flushright}


\section{The space-time RNL-LS parallel algorithm}
\begin{figure}[ht!]
	\centering
\includegraphics[width=1.0\textwidth]{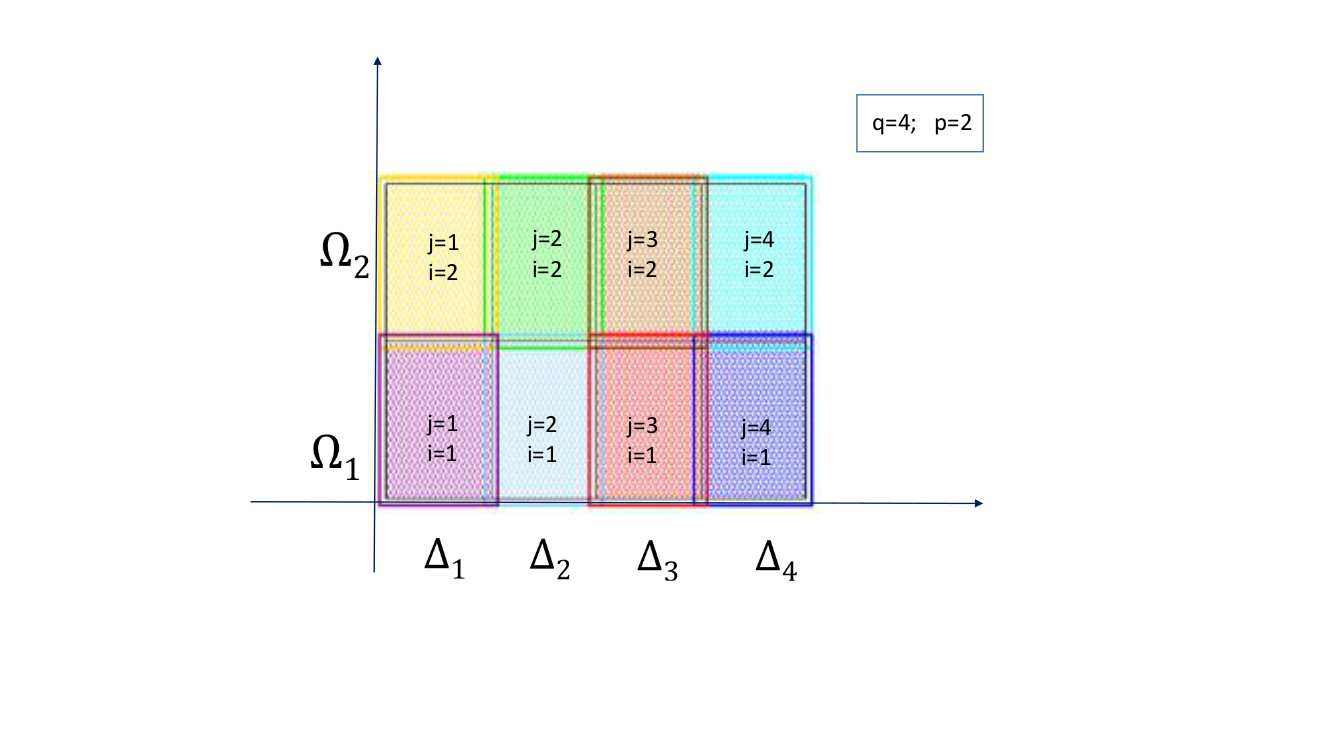}
\caption{Configurations of the decomposition  $DD(\Delta_M \times \Omega_{K})$, if $\Omega \subset \Re$ and $Q=4$, $P=2$.}
\label{imm-conf}
\end{figure}

We  introduce the  algorithm  solving the RNL-LS problem by using the space-time decomposition, in other words,   solving the $QP=q\times p$ local  problems in $\Delta_j \times \Omega_i$, where $j=1,Q$ and $i=1,P$ (see Figure \ref{imm-conf} for an example of domain decomposition where $Q=4$ and $P=2$.).

\begin{definition}[DD-RNL-LS Algorithm]
Let $\mathcal{A}^{loc}_{RNLLS}(\Delta_j \times \Omega_i)$ denote the algorithm solving the local 4D-Var DA problem  defined in $\Delta_j \times \Omega_i$. The  space-time DD-RNL-LS parallel algorithm  solving the RNL-LS problem in $DD(\Delta \times \Omega)$ is symbolically denoted as $$\mathcal{A}^{DD}_{RNNLS}(\Delta_M \times \Omega_{K})$$ and  is defined as the merging  of the $QP=Q\times P$ local algorithms $\mathcal{A}^{loc}_{RNLLS}(\Delta_j \times \Omega_i)$:

\begin{equation}\label{D-TR:algorithm}
\mathcal{A}^{DD}_{RNLLS}(\Delta_M \times \Omega_{NP}):=\bigcup_{j=1,Q;i=1,P} \mathcal{A}^{loc}_{RNLLS}(\Delta_j \times \Omega_i).
\end{equation}
\begin{flushright}
$\spadesuit$
\end{flushright}
\end{definition}

\noindent The DD-RNL-LS algorithm can be sketched as described by \textbf{Algorithm 1}.  Similarly, the Local RNL-LS algorithm $\mathcal{A}^{loc}_{RNLLS}$ is described by \textbf{Algorithm 2}.

\begin{algorithm}{Algorithm 1; $\mathcal{A}^{DD}_{RNLLS}$: solves the RNL-LS problem on $\Delta_M \times \Omega_{NP}$}\label{DD-TRAlg}
\begin{algorithmic}[1]
\Procedure{DD-4DVar}{$in: \mathbf{M}^{\Delta_M \times \Omega_{K}}, \mathbf{u}^b_0, \left\{\mathbf{R}_k\right\}_{k}, \mathbf{B}, \mathbf{H}, \mathbf{v}, \Delta_M, \Omega_{K}; out: \mathbf{u^{\mathbf{DA}}}$}
\State {\% Domain Decomposition Step}
\State Compute $\mathbf{M}^{\Delta_M \times \Omega_{K}}$ from  $\mathcal{M}^{\Delta\times \Omega}  $ 
\State {\% Run  $\mathbf{M}^{\Delta_M \times \Omega_{K}}$} in (\ref{cond1}) with initial condition $\mathbf{u}_0^b$ 
\State $\mathbf{u}^b= \mathbf{M}^{\Delta_M \times \Omega_{K}}(\mathbf{u}^b_0)$
\State {\% Local Model Linearization Step }
\For{$j=1,q;\, i=1, p$}
\State $l:=0, \mathbf{u}_{ji}^{0}= \mathbf{u}_{ji}^{b}$
\Repeat
\State \hspace*{.1cm} $l:=l+1$
\State \hspace*{.2cm}{\bf Call } Loc\_RNLLS $(in: \mathbf{M^{0,M-1}}, \left\{\mathbf{R}_k\right\}_{k}, \mathbf{B}, \mathbf{H}, \mathbf{v}, \mathbf{u}^b, \Delta_j, \Omega_{i}; out: \mathbf{u}^l_{ji})$
\State \hspace*{.2cm}\textbf{Exchange} $\mathbf{u}_{ji}^k$ between adjacent subdomains
\Until{$\|\mathbf{u}_{ji}^l - \mathbf{u}_{ji}^{l-1}\|< eps$}
\EndFor
\State {\% End the Domain Decomposition Step}
\State {\bf Gather} of $\mathbf{u}_{ji}^l: \mathbf{u}^{\mathbf{DA}}= 
{\arg\min}_{ji} \left\{J \left (\mathbf{u}_{ji}^l\right ) \right \}$
 \label{puntoparall2}
 \EndProcedure
\end{algorithmic}
\end{algorithm}

\normalsize
\vspace*{.5cm}


\begin{remark}
We observe that  the $\mathcal{A}^{DD}_{RNNLS}(\Delta_M \times \Omega_{K})$ algorithm is based on  two main steps:  the domain decomposition  step (see line 2) and  the model linearization step (see line 6).  Thus, this algorithm uses a convex approximation of  the objective DA functional so that Theorem  \ref{theorem215} holds.  
\end{remark}

\noindent The common approach for solving RNL-LS problems involves defining a sequence of local approximations of $\mathbf{J}_{ij}$ where each member of the sequence is minimized by  employing Newton's  method or one its variants (such as Gauss--Newton, L-BFGS, or Levenberg--Marquardt). Approximations of $\mathbf{J}_{ij}$ are obtained by expanding  $\mathbf{J}_{ij}$ in a truncated Taylor series, while the  minimum  is obtained by using second-order sufficient conditions \cite{Dennis96,Nocedal}. Let us consider \textbf{Algorithm 2} solving the RNL-LS problem on $\Delta_j \times \Omega_i$.\\[.5cm]

\begin{algorithm}{Algorithm 2; $\mathcal{A}^{loc}_{RNLLS}$: solves an RNL-LS problem  on $\Delta_j \times \Omega_i$ }\label{local-TRAlg:RNL-LS2}
\begin{algorithmic}[1]
\Procedure{Loc-RNLLS}{$in: \mathbf{M^{0,M-1}}, \left\{\mathbf{R}_k\right\}_{k}, \mathbf{B}, \mathbf{H}, \mathbf{v}, \mathbf{u}^b,\Delta_j, \Omega_{i}; out: \mathbf{u}^l_{ji}$}
\State {\bf Initialize} $\mathbf{{u}}^0_{ij}:=\mathbf{{u}}_{ij}^{b}$;
\State {\bf Initialize} $l:=0$;
\State \textbf{repeat} \% at each step  $l$, a  local approximation of $\mathbf{\tilde{J}}_{ij}$  is minimized
\State \hspace*{1cm}{\bf Compute} $\delta \mathbf{{u}}_{ij}^l=\arg\min \,\mathbf{\tilde{J}}_{{ji}}$
\State \hspace*{1cm}{\bf Update} $\mathbf{{u}}_{ji}^l=\mathbf{{u}}_{ji}^l+\delta \mathbf{{u}}_{ji}^l $
\State \hspace*{1cm}{\bf Update} $l=l+1$
\State \textbf{until} (convergence is reached)
 \EndProcedure
\end{algorithmic}
\end{algorithm}

\normalsize
\vspace*{1cm}
\noindent The main computational task occurs at step 5 of \textbf{Algorithm 2} concerning the minimization of $\tilde{\mathbf{J}_{ji}}$, which is the local approximation of $\mathbf{J}_{ij}$.   Two approaches could be employed in \textbf{Algorithm 2}:
\begin{enumerate}
  \item[(a)] By truncating  the Taylor series expansion of $\mathbf{J}_{ij}$ at the second order, we get
    \begin{equation}\label{model}
       \mathbf{J}_{ij}^{QD}(\mathbf{{u}}_{ji}^{l+1})=\mathbf{J}_{ij}(\mathbf{{u}}_{ji}^{l})+\nabla
        \mathbf{J}_{ij}(\mathbf{{u}}_{ji}^{l})^T 
        \delta\mathbf{{u}}_{ji}^l+\left(\delta\mathbf{{u}}_{ji}^l\right)^T \nabla^2 \mathbf{J}_{ij}(\mathbf{{u}}_{ji}^{l}) \delta \mathbf{{u}}_{ji}^l 
    \end{equation}
    giving a quadratic approximation of $\mathbf{J}_{ji}$ at $\mathbf{u}_{ji}^l$. Newton-based methods (including LBFGS and Levenberg--Marquardt) use  $\mathbf{\tilde{J}}_{{ji}}=\mathbf{J}_{ij}^{QD}$.
  \item[(b)]  By truncating the Taylor series expansion of $\mathbf{J}_{ij}$ at the first order, we get the following linear  approximation of $\mathbf{J}_{ij}$ at $\mathbf{{u}}_{ji}^{k}$:
    \begin{equation}\label{modelGN}
      \mathbf{J}_{ij}^{TL}(\mathbf{{u}}_{ji}^{l+1})=\mathbf{J}_{ij}(\mathbf{{u}}_{ji}^{l})+\nabla \mathbf{J}_{ij}(\mathbf{{u}}_{ji}^{l})^T \delta \mathbf{{u}}_{ji}^l=\frac{1}{2} \|\nabla \mathbf{F}_{ji}(\mathbf{{u}}_{ji}^{l})\delta \mathbf{{u}}_{ji}^l+\mathbf{F}_{ji}(\mathbf{{u}}_{ji}^{l})\|_2^2,
    \end{equation}
    where we let\footnote{If $\mathbf{{C}}_{ji}=diag((\mathbf{B}^{-1})_{ji},(\mathbf{R}^{-1})_{ji})$, and $\mathbf{\tilde{d}}_{ji}^{l}=(\mathbf{{u}}_{ji}^{l}-\mathbf{{u}}_{0}^{b},\mathbf{H}_{ji}^{0}(\mathbf{{u}}_{ji}^{l})-\mathbf{{v}}_{ji}^{k}, \ldots, (\mathbf{{H}}^{M-1})_{ji}[(\mathbf{{M}}_{M-2,M-1}^{k})_{ji}(\mathbf{{u}}_{ji}^{k})]-\mathbf{{v}}_{ji}^{l})$, then $\mathbf{J}_{ij}:=\frac{1}{2}((\mathbf{{C}}^{-1/2})_{ji}\mathbf{\tilde{d}}_{ji}^{l})^T((\mathbf{{C}}^{-1/2})_{ji}\mathbf{\tilde{d}}_{ji}^{l}) = \|\mathbf{F}_{ji}\|_2^2$, where $\mathbf{F}_{ji}=(\mathbf{{C}}^{-1/2})_{ji}\mathbf{\tilde{d}}_{ji}^{l}$.
     $\mathbf{J}_{ij}:= \|\mathbf{F}_{ji}\|_2^2$,},
 which gives a linear approximation of $\mathbf{J}_{ji}$ at $\mathbf{u}_{ji}^l$. Gauss--Newton's methods (including truncated or approximated Gauss--Newton \cite{Gratton}) use $\mathbf{J}^{TL}_{ji}=\tilde{\mathbf{J}}_{ji}$.
\end{enumerate}

\noindent Observe that from (\ref{model}) it follows that
    \begin{equation}\label{modelGN2}
     \mathbf{J}_{ij}^{QD}(\mathbf{{u}}_{ji}^{l+1})=\mathbf{J}_{ij}^{
    TL}(\mathbf{{u}}_{ji}^{l})+\frac{1}{2} \left(\delta \mathbf{{u}}_{ji}^l\right)^T  \nabla^2
  \mathbf{J}_{ij}(\mathbf{{u}}_{ji}^{l}) \delta \mathbf{{u}}_{ji}^l.
    \end{equation}
\noindent \textbf{Algorithm 2} can be updated to   \textbf{Algorithm 3} as described below.  \\

\begin{algorithm}{Algorithm 3; $\mathcal{A}^{loc}_{RNLLS}$: solves an RNL-LS problem  on $\Delta_j \times \Omega_i$ }\label{local-TRAlg:RNL-LS3}
\begin{algorithmic}[1]
\Procedure{Loc-RNLLS}{$in: \mathbf{M^{0,M-1}}, \left\{\mathbf{R}_k\right\}_{k}, \mathbf{B}, \mathbf{H}, \mathbf{v}, \mathbf{u}^b,\Delta_M, \Omega_{K}; out: \mathbf{u}^l_{ji}$}
\State {\bf Initialize} $\mathbf{{u}}^0_{ij}:=\mathbf{{u}}_{ij}^{b}$;
\State {\bf Initialize} $l:=0$;
\State \textbf{repeat}
\State \hspace*{1cm}{\% Compute} $\delta \mathbf{{u}}_{ij}^l=\arg\min \,\mathbf{J}_{{ji}}$ by using $\mathcal{A}_{QN}^{loc}$ or $\mathcal{A}_{LLS}^{loc}$
\State \hspace*{1cm}  \textbf{If} (QN) \textbf{then}
\State \hspace*{1.5cm}{\bf Call} Loc-QN ({$in: \mathbf{M^{0,M-1}}, \left\{\mathbf{R}_k\right\}_{k}, \mathbf{B}, \mathbf{H}, \mathbf{v}, \mathbf{u}^b,\Delta_M, \Omega_{K}; out: \mathbf{u}^l_{ji}$ })
\State \hspace*{1cm} \textbf{ElseIf} (LLS) \textbf{then}
\State \hspace*{1.5cm}{\bf Call} Loc-LLS ({$in: \mathbf{M^{0,M-1}}, \left\{\mathbf{R}_k\right\}_{k}, \mathbf{B}, \mathbf{H}, \mathbf{v}, \mathbf{u}^b,\Delta_{M}, \Omega_{K}; out: \mathbf{u}^l_{ji}$})
\State \hspace*{1cm}  \textbf{EndIf}
\State \hspace*{1cm}{\bf Update} $\mathbf{{u}}_{ji}^l=\mathbf{{u}}_{ji}^l+\delta \mathbf{{u}}_{ji}^l $
\State \hspace*{1cm}{\bf Update} $l=l+1$
\State \textbf{until} (convergence is reached)
 \EndProcedure
\end{algorithmic}
\end{algorithm}

 \begin{enumerate}
   \item[(a)] $\mathcal{A}_{QN}^{loc}$: computes a local minimum of $\mathbf{J}^{QN}_{ji}$ following  the Newton  descent direction. The minimum is computed by solving the linear system involving the Hessian matrix $\nabla ^2\mathbf{{J}}_{ij}$ and the negative gradient $-\nabla \mathbf{{J}}_{ij}$  at $\mathbf{{u}}_{ji}^l$, for each value of $l$ (see \textbf{Algorithm 4} described below).\\

\begin{algorithm}{Algorithm 4; $\mathcal{A}^{loc}_{QLS}$: solves a Q-LS problem  on $\Delta_j \times \Omega_i$}\label{local-TRAlg:Q-LS}
\begin{algorithmic}[1]
\Procedure{Loc-QN}{$\mathbf{M^{0,M-1}}, \left\{\mathbf{R}_k\right\}_{k}, \mathbf{B}, \mathbf{H}, \mathbf{v}, \mathbf{u}^b, \Delta_M, \Omega_{K}; out: \mathbf{u}^l_{ji}$}
\State {\bf Initialize} $\mathbf{{u}}^0_{ji}:=\mathbf{{u}}_{ji}^{b}$;
\State {\bf Initialize} $l:=0$;
\State \textbf{repeat}
\State \hspace*{.51cm}{\%Compute} $\delta \mathbf{{u}}_{ij}^l=\arg\min \,\mathbf{J}_{{ji}}^{QD}$, by Newton's method
\State \hspace*{.51cm}{\bf 1.1 Compute} $\nabla \mathbf{J}_{{ji}}(\mathbf{{u}}_{ij}^l)= \nabla \mathbf{F}_{{ji}}^T(\mathbf{{u}}_{ji}^l)\nabla \mathbf{F}_{{ji}}(\mathbf{{u}}_{ji}^l)$
\State \hspace*{.51cm}{\bf 1.2 Compute} $\nabla^2 \mathbf{J}_{{ji}}(\mathbf{{u}}_{ij}^l)=\nabla \mathbf{F}_{{ji}}^T(\mathbf{{u}}_{ji}^l)\nabla \mathbf{F}_{{ji}}(\mathbf{{u}}_{ji}^l)+\mathbf{Q}((\mathbf{{u}}_{ij}^l))$
\State \hspace*{.51cm}{\bf 1.3 Solve} $\nabla^2 \mathbf{J}_{{ji}}(\mathbf{{u}}_{ij}^l) \delta \mathbf{{u}}_{ij}^l=-\nabla \mathbf{J}_{{ji}}(\mathbf{{u}}_{ij}^l)$
\State \hspace*{.51cm}{\bf Update} $\mathbf{{u}}_{ji}^l=\mathbf{{u}}_{ji}^l+\delta \mathbf{{u}}_{ji}^l $
\State \hspace*{.51cm}{\bf Update} $l=l+1$
\State \textbf{until} (convergence is reached)
 \EndProcedure
\end{algorithmic}
\end{algorithm}

\item[(b)]$\mathcal{A}_{LLS}^{loc}$: computes a local minimum of $\mathbf{J}^{TL}_{ji}$ following the steepest descent direction. The minimum is computed by solving the normal equations arising from  the local linear least squares (LLS) problem (see \textbf{Algorithm 5} described below).\\[.5cm]

\begin{algorithm}{Algorithm 5; $\mathcal{A}^{loc}_{LLS}$: solves LLS problems in $\Delta_j \times \Omega_i$  }\label{local-TRAlg:LLS}
\begin{algorithmic}[1]
\Procedure{Loc-LLS}{$\mathbf{M^{0,M-1}}, \left\{\mathbf{R}_k\right\}_{k}, \mathbf{B}, \mathbf{H}, \mathbf{v}, \mathbf{u}^b, \Delta_M, \Omega_{K}; out: \mathbf{u}^l_{ji}$}
\State {\bf Initialize} $\mathbf{{u}}^0_{ij}:=\mathbf{{u}}_{ij}^{b}$;
\State {\bf Initialize} $l:=0$;
\State \textbf{repeat}
\State \hspace*{.5cm}{\bf Compute} $\nabla \mathbf{J}_{{ji}}=\nabla \mathbf{F}_{{ji}}^T(\mathbf{{u}}_{ji}^l)\nabla \mathbf{F}_{{ji}}(\mathbf{{u}}_{ji}^l)$
\State \hspace*{.5cm}{\%Compute} $\delta \mathbf{{u}}_{ij}^l=\arg\min \,\mathbf{J}_{{ji}}^{TL}$ by solving the normal equations system:
\State \hspace*{.5cm}{\bf Solve} $\nabla \mathbf{F}_{{ji}}^T(\mathbf{{u}}_{ji}^l)\nabla \mathbf{F}_{{ji}}(\mathbf{{u}}_{ji})\delta \mathbf{{u}}_{ji}^l=-\nabla \mathbf{F}_{{ji}}^T(\mathbf{{u}}_{ji}^l)\mathbf{F}_{{ji}}(\mathbf{{u}}_{ji}^l)$
\State \hspace*{.5cm}{\bf Update} $\mathbf{{u}}_{ji}^l=\mathbf{{u}}_{ji}^l+\delta \mathbf{{u}}_{ji}^l $
\State \hspace*{.5cm}{\bf Update} $l=l+1$
\State \textbf{until} (convergence is reached)
 \EndProcedure
\end{algorithmic}
\end{algorithm}
\end{enumerate}

\begin{remark}: We observe that if, in the  $\mathcal{A}^{loc}_{QN}$ algorithm, matrix $\mathbf{Q}(\mathbf{{u}}_{ij}^l)$  (see line 6 of Algorithm 4) is neglected, we get the Gauss--Newton method described by $\mathcal{A}^{loc}_{LLS}$ algorithm. More generally,
the term $\mathbf{Q}(\mathbf{{u}}_{ij}^l)$ 
\begin{enumerate}
  \item in the case of Gauss--Newton, $\mathbf{Q}(\mathbf{{u}}_{ij}^l)$, is neglected;
  \item in the case of Levenberg--Marquardt, $\mathbf{Q}(\mathbf{{u}}_{ij}^l)$ equals  $\lambda I$, where the damping term,  $\lambda >0$, is updated at each iteration and $I$ is the identity matrix \cite{Lev,Marq}; and
  \item in the case of the L-BFGS, the Hessian matrix is rank-1 updated at every iteration \cite{LBFGS}. 
\end{enumerate}
\end{remark}
\noindent 

\noindent In accordance with the most common implementation of the 4D-Var DA \cite{Dati,ROMS}, we focus  attention on the Gauss--Newton(G-N) method described in $\mathcal{A}^{loc}_{LLS}$ in \textbf{Algorithm 6}.\\

  \noindent For each $l$, let  $\mathbf{G}_{ji}^l = RO_{ji}[\mathbf{G}^l]$, where $\mathbf{G}^l \in \Re^{( M\times nobs)\times( NP \times M)}$, be the block diagonal matrix such that
  \begin{equation}\label{G}
    \mathbf{G}^l= \left \{
    \begin{array}{ll}
     diag\, [\mathbf{H}^0,\mathbf{H}^1 \mathbf{M}^{0,1}_l, \ldots, \mathbf{H}^{M-1} \mathbf{M}^{M-2,M-1}_l]  & M>1; \\
       \mathbf{H}^0& M=1,
    \end{array}
    \right .
  \end{equation}
  where $(\mathbf{G}_{ji}^T)^l=RO_{ji}[(\mathbf{G}^T)^l]$ is the restriction of the transpose of $\mathbf{G}^l$  and $$\mathbf{M}^{0,1}_l, \ldots, \mathbf{M}^{M-2,M-1}_l$$ are the TLMs of $\mathbf{M}^{k,k+1}$, for $s=0,M-1$, around $ \mathbf{u}^l_{ji}$, respectively. 
  Let  $$\mathbf{d}_{ji}^l=\mathbf{v}_{ji}-\mathbf{H}_{ji}\mathbf{u}^l_{ji}\quad $$
be the restriction of the misfit vector where $\mathbf{H}_{ji}$ is the matrix 

\begin{displaymath}
\mathbf{H}_{ji}=diag\, [ (RO_{ji}[\mathbf{H}_{ji}^{k}])_{k: t_k \in\Delta j} ].
\end{displaymath}

\noindent Let  $\mathbf{R}_{ji}$ the block diagonal matrix such that
\begin{displaymath}
\mathbf{R}_{ji}=diag\, [(RO_{ji}[\mathbf{R}^k ])_{k: t_k \in\Delta j} ].
\end{displaymath}

\noindent In line 7 of \textbf{Algorithm 5},  it is
\begin{equation}\label{line5and7}
 \nabla \mathbf{F}_{{ji}}^T(\mathbf{{u}}_{ji}^l)\nabla \mathbf{F}_{{ji}}(\mathbf{{u}}_{ji}^l)=\mathbf{B}^{-1}_{ji}+(\mathbf{G}_{ji}^T)^l\mathbf{R}_{ji}\mathbf{G}_{ji}^l,
\end{equation}

\noindent and

\begin{equation}\label{line6and7}
  -\nabla \mathbf{F}_{{ji}}^T(\mathbf{{u}}_{ji}^l) \mathbf{F}_{{ji}}(\mathbf{{u}}_{ji}^l)=(\mathbf{G}_{ji}^T)^l\mathbf{R}_{ji}^{-1}\mathbf{d^l}_{ji},
\end{equation}

\noindent where $B_{ji}$ and $R_{ji}$ are the restrictions of $B$ and $R$ matrices, respectively.

  \noindent  
 Most popular 4D-Var DA software implements the so-called $\mathbf{B}$-preconditioned Krylov subspace iterative method \cite{Gratton,Gurol,ROMS} arising  by using the background error covariance matrix as a preconditioner of a Krylov subspace iterative method. \\

 \noindent Let ${\bf B}_{ji}=\mathbf{V}_{ji}\mathbf{V}_{ji}^T$
 be expressed in terms of the deviance matrix $\mathbf{V}_{ji}$ and $\mathbf{w}_i$ such that
\begin{equation}\label{EqV}
 \mathbf{w}_{ji}^l =  \mathbf{V}_{ji} ^+ (\mathbf{u}_{ji}^l-\mathbf{u}_{ji}^b)
\end{equation}
with $V_i ^+$ the generalized inverse of $\mathbf{V}_i $. Then (\ref{line5and7}) becomes
\begin{equation}\label{Eq_matrixAV:1}
\mathbf{B}^{-1}_{ji}+(\mathbf{G}_{ji}^T)^l\mathbf{R}_{ji}\mathbf{G}_{ji}^l={\bf I}_{ji}+({\bf G}_{ji}^l\mathbf{V}_{ji})^T({\bf R}^{-1})_{ji}{\bf G}_{ji}^l\mathbf{V}_{ji},
\end{equation}
and(\ref{line6and7}) becomes
\begin{equation}\label{Eq_matrixAV:2}
(\mathbf{G}_{ji}^T)^l(\mathbf{R}^{-1})_{ji}\mathbf{d}_{ji}=({\bf G}_{ji}\mathbf{V}_{ji})^T)^k({\bf R}^{-1})_{ji}\mathbf{d}_{ji}.
\end{equation}
The  normal equation system (see line 7 of $\mathcal{A}^{loc}_{LLS}$), in other words, the linear system
$$((\mathbf{B}^{-1})_{ji}+(\mathbf{G}_{ji}^T)^l\mathbf{R}_{ji}\mathbf{G}_{ji}^l  )\delta \mathbf{{u}}_{ji}^l=(\mathbf{G}_{ji}^T)^l(\mathbf{R}^{-1})_{ji}\mathbf{d}_{ji},$$ 
becomes
$$({\bf I}_{ji}+({\bf G}_{ji}^l\mathbf{V}_{ji})^T({\bf R}^{-1})_{ji}{\bf G}_{ji}^l\mathbf{V}_{ji})\delta \mathbf{{u}}_{ji}^l=({\bf G}_{ji}^l\mathbf{V}_{ji})^T({\bf R}^{-1})_{ji}\mathbf{d}_{ji}\quad .$$

\begin{definition}[DD-4D-Var Algorithm]
Let $\mathcal{A}^{loc}_{4DVar}(\Delta_j \times \Omega_i)$ denote the algorithm solving the local 4D-Var DA problem  defined in $\Delta_j \times \Omega_i$.  The space-time 4D-Var DA parallel algorithm  solving the 4D-Var DA problem in $DD(\Delta_{M} \times \Omega_{K})$ is symbolically denoted as $\mathcal{A}^{DD}_{4DVar}(\Delta_M \times \Omega_{K})$, and it is defined as the union of the $QP=q\times p$ local algorithms $\mathcal{A}^{loc}_{4DVar}(\Delta_j \times \Omega_i)$:

\begin{equation}\label{D-TR:algorithm:2}
\mathcal{A}^{DD}_{4DVar}(\Delta_M \times \Omega_{K}):=\bigcup_{j=1,q;i=1,p} \mathcal{A}^{loc}_{4DVar}(\Delta_j \times \Omega_i).
\end{equation}
\begin{flushright}
$\spadesuit$
\end{flushright}
\end{definition}

\noindent Algorithm   $ \mathcal{A}^{loc}_{4DVar}$ is Algorithm $ \mathcal{A}^{loc}_{LLS}$ (see \textbf{Algorithm 5}) specialized for the 4D-Var DA problem, and it  is  described by \textbf{Algorithm 6} and \textbf{Algorithm 7}, described below \cite{Gurol}.\\

\begin{algorithm}{Algorithm 6; $\mathcal{A}^{loc}_{4DVar}$: solves  Local 4DVAR DA problem in $\Delta_j \times \Omega_i$  }\label{local-TRAlg:4DVar}
\begin{algorithmic}[1]
\Procedure{Loc-4DVar}{$\mathbf{M}^{\Delta_M\times \Omega_{NP}}, \mathbf{R}, \mathbf{B}, \mathbf{H}, \mathbf{v}, \mathbf{u}^b, \Delta_M, \Omega_{K}; out: \mathbf{u}^l_{ji}$}
\State {\bf Initialize} $\mathbf{{u}}^0_{ji}:=\mathbf{{u}}_{ji}^{b}$;
\State {\bf Initialize} $l:=0$;
\State \textbf{repeat}
\State \hspace*{.5cm}{\bf Compute} $\mathbf{{d}}_{ji}^l=\mathbf{{v}}_{ji}-\mathbf{H}_{ji}(\mathbf{{u}}_{ji}^l)$
\State \hspace*{.5cm}\textbf{Call} TLM($in: \,\mathcal{M}^{\Delta \times \Omega}$, $\mathbf{u}_{ji}^{l}$; $out:\mathbf{M}_{0,M-1}^l$)
\State \hspace*{.5cm}\textbf{Call} ADJ($in:\,\mathbf{M}^k_{0,M-1}$; $out:(\mathbf{M}_{0,M-1}^T)^l$)
\State \hspace*{.5cm}{\bf Compute} ${\bf G}_{ji}$, $\mathbf{V}_{ji}$
\State \hspace*{.5cm}{\bf Call} $ \mathcal{A}^{loc}_{BLanczos}$ ($\mathbf{G}_{ji}^k, \mathbf{V}_{ji}, \mathbf{R}_{ji}, \mathbf{B}_{ji}, \mathbf{d}_{ji}, \mathbf{u}^b_{ji}, \Delta_j, \Omega_i; out: \delta\mathbf{u}^k_{ji}$)
\State \hspace*{.5cm}{\bf Update} $\mathbf{{u}}_{ji}^l=\mathbf{{u}}_{ji}^l+\delta \mathbf{{u}}_{ji}^l $
\State \hspace*{.5cm}{\bf Update} $l=l+1$
\State \textbf{until} (convergence is reached)
 \EndProcedure
 \State \textbf{endprocedure}
 \Procedure{TLM}{$in: \,\mathcal{M}^{\Delta \times \Omega}$, $\mathbf{u}_{ji}^{l}$; $out:\mathbf{M}_{0,M-1}^l$}
\State \hspace*{.5cm}{\%Linearize $\mathbf{M}^{\Delta_M \times \Omega_{NP}}$ about $\mathbf{u}_{ji}^{l}$}
\EndProcedure
\State \textbf{endprocedure}
\Procedure{ADJ}{$in:\,\mathbf{M}^k_{0,M-1}$; $out:(\mathbf{M}_{0,M-1}^T)^l$}
\State \hspace*{.5cm}{\%Compute the adjoint of  $\mathbf{M}_{0,M-1}$}
 \EndProcedure
 \State \textbf{endprocedure}
\end{algorithmic}
\end{algorithm}

\begin{algorithm}{Algorithm 7; $\mathcal{A}^{loc}_{BLanczos}$: \textbf{B}Lanczos for 4D-VAR DA problem in $\Delta_j \times \Omega_i$ }\label{local-TRAlg:BLanczos}
\begin{algorithmic}[1]
\Procedure{BLanczos-4DVar}{$\mathbf{G}_{ji}, \mathbf{V}_{ji}, \mathbf{R}_{ji}, \mathbf{B}_{ji}, \mathbf{d}_{ji}, \mathbf{u}^b_{ji}, \Delta_j, \Omega_i; out: \delta\mathbf{u}^l_{ji}$}
\State \hspace*{.5cm}{\% Solve} $ ({\bf I}_{ji}+({\bf G}_{ji}\mathbf{V}_{ji})^T({\bf R}^{-1})_{ji}{\bf G}_{ji}\mathbf{V}_{ji})\delta \mathbf{{u}}_{ji}^l=({\bf G}_{ji}\mathbf{V}_{ji})^T({\bf R}^{-1})_{ji}\mathbf{d}_{ji}$
\State \hspace*{.5cm}{\%    by using BLanczos algorithm (see \cite{Gurol}) }
 \EndProcedure
\end{algorithmic}
\end{algorithm}
\noindent In the next section we will show that this formulation leads to local numerical solutions that converge to the numerical solution of the global problem. 

\section{Convergence analysis} 
In the following we assume $\|\cdot\|=\|\cdot \|_{\infty}$.

\begin{proposition}\label{convergenza}
Let $u_{j,i}^{ASM,r}$ be the approximation of the increment $\delta \mathbf{{u}}_{ji}$ to the solution $\mathbf{{u}}_{ji}$  obtained at step $r$ of ASM-based inner loop on $\Omega_{j}\times \Delta_{i}$. Let $u_{j,i}^{l}$ be the approximation of $\mathbf{u}_{j,i}$ obtained at step $l$ of the outer loop, that is,   the space-time decomposition approach  on $\Omega_{j}\times \Delta_{i}$.  Let us assume that 
 the numerical scheme   discretizing  the model $\mathbf{M}_i^{j,j+1}$  is convergent.
Then with fixed  $i$ and $j$, it holds that
\begin{equation}\label{tesi0}
\forall \epsilon>0\ \ \exists M(\epsilon)>0 \ \ : \ \ l>M(\epsilon) \ \ \Rightarrow \ \ 
E_{j,i}^{l}:=\|\mathbf{u}_{j,i}-u_{j,i}^{l}\| \le \epsilon.
\end{equation}

\end{proposition}

\noindent \textbf{Proof}: Let $u_{j,i}^{\mathbf{M}_i^{j,j+1},l+1}$ be the numerical solution of $\mathbf{M}_i^{j,j+1}$ at step $l$; taking into account that, according to the incremental update of the solution of  the 4D-Var DA functional (for instance, see line 10 of Algorithm 7), the approximation $\mathbf{u}_{j,i}^l$  is computed as 
$$\mathbf{u}_{j,i}^l= u_{j,i}^{\mathbf{M}_i^{j,j+1},l+1}+[u_{j,i}^{ASM,r}-u_{j,i}^{\mathbf{M}_i^{j,j+1},l}] , $$
and then  
\begin{equation}\label{pri}
\begin{array}{ll}
E_{j,i}^{l}:=\|\mathbf{u}_{j,i}-u_{j,i}^{l}\|&=\|\mathbf{u}_{j,i}-u_{j,i}^{\mathbf{M}_i^{j,j+1},l+1}-[u_{j,i}^{ASM,r}-u_{j,i}^{\mathbf{M}_i^{j,j+1},l}]\|\\&\le \|\mathbf{u}_{j,i}-u_{j,i}^{ASM,r}\|+\|{u}_{j,i}^{\mathbf{M}_i^{j,j+1},l}-u_{j,i}^{\mathbf{M}_i^{j,j+1},l+1} \|\,.
\end{array}
\end{equation}
From the hypothesis above we have 
   \begin{eqnarray}
\forall \epsilon^{\mathbf{M}_i^{j,j+1}}>0 \,,~ \exists \,M^{1}(\epsilon^{\mathbf{M}_i^{j,j+1}})>0 & : &  l>M^{1}(\epsilon^{\mathbf{M}_i^{j,j+1}}) \nonumber\\
& \Rightarrow &  \|u_{j,i}^{\mathbf{M}_i^{j,j+1},l+1}-u_{j,i}^{\mathbf{M}_i^{j,j+1},l}\| \le \epsilon^{\mathbf{M}_i^{j,j+1}}, \label{ipotesi1}
\end{eqnarray}
and (\ref{pri}) can be rewritten as follows:
\begin{equation}\label{rel}
\begin{array}{ll}
\|\mathbf{u}_{j,i}-u_{j,i}^{l}\|\le \|\mathbf{u}_{j,i}-u_{j,i}^{ASM,r}\|+\epsilon^{\mathbf{M}_i^{j,j+1}}.
\end{array}
\end{equation}
Convergence of ASM is proved in \cite{Clerc}. Similarly, applying ASM to the 4D-Var DA problem,  we have that 
\begin{equation}\label{tesi}
\forall \epsilon^{ASM}>0\ \ \exists M^{2}(\epsilon^{ASM})>0 \ \ : \ \ r >M^{2}(\epsilon^{ASM}) \ \ \Rightarrow \ \ 
\|u_{j,i}-u_{j,i}^{ASM,r}\| \le \epsilon^{ASM},
\end{equation}
and for $l\ >M^{2}(\epsilon^{ASM})$, we get 
\begin{equation}\label{rel1}
\begin{array}{ll}
\|u_{j,i}-u_{j,i}^{ASM,l}\|\le \epsilon^{ASM} +\epsilon^{\mathbf{M}_i^{j,j+1}}.
\end{array}
\end{equation}
Hence, by using $\epsilon:={\epsilon^{ASM}}+\epsilon^{\mathbf{M}_i^{j,j+1}}$ and $M(\epsilon):= \max\{M^{1}(\epsilon^{ASM}),M^{2}(\epsilon^{\mathcal{M}_i^{j,j+1}})\}$, we obtain (\ref{tesi}).\\

\normalsize
\noindent Convergence behavior of local solutions essentially  depends on the rate of convergence of the truncation error given by the discrete forecasting model (see \cite{report_conv} for the convergence analysis).

\section{Performance Analysis}
  We  use time complexity and scalability as performance metrics. Our aim is to highlight the benefits arising from using the  decomposition approach instead of solving the problem on the whole domain. As we  discuss later, the performance gain that we get from using the space and time decomposition  approach is twofold.
\begin{enumerate}
  \item Instead of solving one larger problem, we can solve several smaller problems that are better conditioned than the former problem.  This approach leads to a reduction in each local algorithm’s time complexity.
  \item Subproblems reproduce the whole problem at smaller dimensions, and they are solved in parallel. This approach leads to a reduction in software execution time.
\end{enumerate}

\noindent We  give  the following definition.
\noindent  \begin{definition}
A uniform bidirectional decomposition of the space and time domain $\Delta_M \times \Omega_{K}$ is such that  if we let $$size(\Delta_M \times \Omega_{K})= M \times K $$ be the size of the whole domain, then each subdomain $\Delta_j \times \Omega_i$ is such that  $$size(\Delta_j \times \Omega_i)=D_t \times D_s, \quad j=1,\ldots,q; \quad i=1,\ldots,p ,$$
where $D_t=\frac{M}{q}\geq 1$  and $D_s=\frac{K}{p}\geq 1$.
\begin{flushright}
$\spadesuit$
\end{flushright} 
 \end{definition}
 In the following we let 
 $$N:= M \times K; \quad N_{loc}:= D_t \times D_s  ; \quad QP:= q \times p \, . $$
\noindent  Let  $T(\mathcal{A}^{DD}_{4DVar}(\Delta_M \times \Omega_{K}))$ denote  time complexity of  $\mathcal{A}^{DD}_{4DVar}(\Delta_{M} \times \Omega_{K})$.  
\noindent We now provide an estimate of the time complexity of each local algorithm, denoted as $T(\mathcal{A}^{Loc}_{4DVar}(\Delta_j \times \Omega_i))$.  This algorithm consists of two loops: an outer loop, over $l$-index, for computing local approximations of $\mathbf{J}_{ji}$, and an inner loop over the $m$ index, for performing the Newton or Lanczos steps. The major computational task to be performed at each step of the outer loop is the computation of $\mathbf{J}_{ji}$.  The major computational tasks to be performed at each step $l$ of the inner loop, in the case of the G-N method (see Algorithm $\mathcal{A}^{Loc}_{4DVar}$), involving the predictive model,  are as follows:\footnote{These assumptions hold true for the so-called local discretization schemes, i.e., those schemes where  each grid point receives contribution from  a neighborhood (for instance,  using finite difference and finite volume discretization schemes as in \cite{Shchepetkin}).}
\begin{enumerate}
  \item Computation of the tangent linear model  $RO_{ji}[\mathbf{M}^{k,k+1}]$ (the time complexity of such an operation scales as the problem size squared) 
  \item Computation of the adjoint model $RO_{ji}[(\mathbf{M}^{k,k+1})^T]$,  which is at least $4$ times more expensive than the computation of $RO_{ji}[\mathbf{M}^{k,k+1}]$
  \item Solution of the normal equations, involving at each iteration two matrix-vector products with $RO_{ji}[(\mathbf{M}^{k,k+1})^T]$   and $RO_{ji}[\mathbf{M}^{k,k+1}]$ (whose time complexity scales as the problem size squared).
\end{enumerate}
\noindent
\noindent Since the most time-consuming operation involving the predictive model is the computation of the tangent linear model,  we prove the following.
\begin{proposition}
Let  $$ P(N_{loc})= a_{d} N_{loc}^{d} + a_{d-1} N_{loc}^{d-1}+ \ldots + a_0, \quad a_d \neq 0$$
 be the polynomial of degree $d=2$ denoting the time complexity of the tangent linear model  $RO_{ji}[\mathbf{M}^{k,k+1}]$.
 Let $m_{ji}$ and  $l_{ji}$ be the number of steps of the outer/inner loop of $\mathcal{A}^{Loc}_{4DVAR}$, respectively.  We get

 $$T(\mathcal{A}^{Loc}_{4DVAR}(\Delta_j \times\Omega_i)))= O\left (m_{ji}l_{ji}P(N_{loc})\right ).$$
 \noindent \textbf{Proof:} It is
 \begin{eqnarray}
  \nonumber
   T(\mathcal{A}^{Loc}_{4DVAR}(\Delta_j \times\Omega_i)) & =  & \\ \nonumber l_{ji} \times \left [T(RO_{ji}[\mathbf{M}^{k,k+1}]) + m_{ji}\times O\left (T(RO_{ji}[\mathbf{M}^{k,k+1}])+T(RO_{ji}[(\mathbf{M}^{k,k+1})^T])\right )\right ] & = & \\ \nonumber
   l_{ji} \times \left [T(RO_{ji}[\mathbf{M}^{k,k+1}]) + m_{ji}\times O\left (T(RO_{ji}[\mathbf{M}^{k,k+1}])+T(RO_{ji}[(\mathbf{M}^{k,k+1})^T])\right )\right ] &=& \\ \nonumber = O\left ( m_{ji}l_{ji}P(N_{loc})\right ).&\\
 \end{eqnarray}

\begin{flushright}
$\clubsuit$
\end{flushright}
\end{proposition}
Let $$m_{max}:= \max_{ji} \,m_{ji};\quad l_{max}:= \max_{ji} \,l_{ji}.$$
Observe that  $m_{max}$ and  $l_{max}$ actually are the number of steps of the outer and inner loops of $\mathcal{A}^{DD}(\Delta_M \times \Omega_{K})$, respectively.  
Let $\mathcal{A}^{G}(\Delta_M \times \Omega_{K})$ denote the algorithm used to solve problem (\ref{min3}) on the undecomposed domain, and 
let $m_G$ and $l_G$  denote the number of iterations of the inner and outer loop of $\mathcal{A}^{G}(\Delta_M \times \Omega_{K})$  algorithm, respectively. Then we have  the following.
\begin{definition} Let 
 $$ \rho^G:= m_G \times l_G\quad ;\quad  \rho^{ji}:= m_{ji} \times l_{ji}\quad ;\quad  \rho^{DD}:= m_{max} \times l_{max}\quad $$
denote the total number of iterations of $\mathcal{A}^{G}_{4DVAR}(\Delta_M \times \Omega_{K})$, of $\mathcal{A}^{Loc}_{4DVAR}(\Delta_j \times \Omega_{i})$  and of $\mathcal{A}^{DD}_{4DVAR}(\Delta_M \times \Omega_{K})$,  respectively.\\[.5cm]
 \end{definition}
If we denote by $\mu(J)$ the condition number of the DA operator, since it holds that \cite{JCP2017} $$\forall \, i,j \quad \mu(J_{4DVAR}^{Loc}) < \mu(J_{4DVAR}),$$
then  $$ \rho^{ji}< \rho^G,$$ and 
$$\rho^{DD} < \rho^G\quad .$$

\noindent This result says that the number of iterations of the $\mathcal{A}^{DD}_{4DVar}(\Delta_M \times \Omega_{K})$ algorithm is always smaller than the number of iterations of the  $\mathcal{A}^{G}_{4DVar}(\Delta_M \times \Omega_{K})$ algorithm. This is one of the benefits of using the space and time decomposition.


\noindent Algorithm scalability  is measured in terms of  \emph{strong scaling} (which is the measure of the algorithm's capability to exploit performance of high-performance computing architectures in order to minimise the time to solution for a given problem with a fixed dimension) and  of  \emph{weak scaling} (which is the measure of the algorithm's capability to use additional computational resources effectively to solve increasingly larger problems).  Various metrics have been developed to
assist in evaluating the scalability  of a parallel algorithm;  speedup, model
throughput, scale-up, efficiency are the most used. Each one highlights  specific needs and limits to be answered by the parallel algorithm. In our case, since we focus mainly on the benefits arising from the use of hybrid computing architectures,  we  consider the so-called  scale-up factor first introduced in  \cite{DD-DA}.\\
\noindent 
\noindent The first result straightforwardly derives from the definition  of the scale-up factor: 
\noindent \begin{proposition}[DD-4D-Var Scale-up factor]
The (relative) scale-up factor of  $\mathcal{A}^{DD}_{4DVar}(\Delta_M \times \Omega_{K})$ related to $\mathcal{A}^{loc}_{4DVar}(\Delta_j \times \Omega_i)$, denoted as
$Sc_{QP}(\mathcal{A}^{DD}_{4DVar}(\Delta_M \times \Omega_{K}))$, is  $$Sc_{QP}(\mathcal{A}^{DD}(\Delta_M \times \Omega_{K})):=\frac{1}{QP}\times\frac{T(\mathcal{A}^{G}_{4DVar}(\Delta_M \times \Omega_{K}))}{ T(\mathcal{A}^{loc}_{4DVar}(\Delta_j \times \Omega_i))}\,\,\,, $$
where $QP:=q \times p$ is the number of subdomains. It is

\begin{equation}\label{scaleup}
  Sc_{QP}(\mathcal{A}^{DD})\geq \frac{\rho^{G}}{\rho^{DD}} \alpha(N_{loc},QP)\,(QP)^{d-1},
\end{equation}
where
$$\alpha(N_{loc},QP)=\frac{a_d+a_{d-1}\frac{1}{N}+\ldots+\frac{a_0}{N_{loc}^{d}}}{a_d+a_{d-1}\frac{QP}{N_{loc}}+\ldots+\frac{a_0(QP)^{d}}{N_{loc}^{d}}}  \quad $$
and
$$ \lim_{QP\rightarrow N_{loc}} \alpha(N_{loc},QP)=\beta \in]0,1].$$
\begin{flushright}
$\spadesuit$
\end{flushright}
\end{proposition}

\begin{corollary}
If $ a_i= 0 \quad \forall i\in[0,d-1]$, then $\beta=1$, that is, 
$$ \lim_{QP\rightarrow N_{loc}} \alpha(N_{loc},QP)=1.$$
Then,
$$\lim_{N_{loc}\rightarrow \infty} \alpha(N_{loc},QP)=1.$$
\begin{flushright}
$\clubsuit$
\end{flushright}

\end{corollary}

\begin{corollary}\label{corscaleup}
If $N_{loc}$ is fixed, then
$$ \lim_{QP\rightarrow N_{loc} }Sc_{1,QP}(\mathcal{A}^{DD})= \beta\cdot N_{loc}^{d-1} \quad ;$$
while if $QP$ is fixed, then
$$ \lim_{N_{loc}\rightarrow \infty }Sc_{1,QP}(\mathcal{A}^{DD})= const \neq 0\quad . $$
\begin{flushright}
$\clubsuit$
\end{flushright}

\end{corollary}

\noindent  From (\ref{scaleup}) it results that, considering one iteration of the whole parallel algorithm,  the growth of the scale-up factor essentially is one order less than the time complexity of the reduced model. In other words, the time complexity of the reduced model impacts mostly the scalability of the parallel algorithm.  In particular, since parameter $d$ is equal to $2$, it follows that the asymptotic scaling factor of the parallel algorithm, with respect to $QP$, is bounded above by two.\\

\noindent Besides the time complexity, scalability is also affected by the  communication overhead of the parallel algorithm. The surface-to-volume ratio is a measure of the amount of data exchange (proportional to surface area of domain) per unit operation (proportional to volume of domain). We prove the following.
\begin{theorem}
 The surface-to-volume ratio of a uniform bidimensional  decomposition of the space-time domain $\Delta_M \times \Omega_{K}$ is
\begin{equation}\label{surfacetovolume}
 \frac{\mathcal{S}}{\mathcal{V}}(\mathcal{A}^{loc}_{4DVar})=2\left (\frac{1}{D_t}+\frac{1}{D_s} \right)\quad .
\end{equation}
\end{theorem}
\noindent
    Let $\mathcal{S}(\mathcal{A}^{loc}_{4DVar})$ denote the surface of each subdomain. Then
     $$\mathcal{S}(\mathcal{A}^{loc}_{4DVar})= 2\left (\frac{M}{q}+\frac{K}{p}\right ) $$
     and  $\mathcal{V}(\mathcal{A}^{loc}_{4DVar})$ denote its volume. Then
     $$\mathcal{V}(\mathcal{A}^{loc}_{4DVar})= \frac{M}{q}\times \frac{K}{p}  \quad .$$
\noindent It holds that
$$\frac{\mathcal{S}}{\mathcal{V}}(\mathcal{A}^{loc}_{4DVar})=\frac{2\left (\frac{M}{q}+ \frac{K}{p}\right )}{ \frac{M}{q}\times \frac{K}{p}}=2\left(\frac{1}{D_t}+\frac{1}{D_s}\right),$$
and (\ref{surfacetovolume}) follows.\\

\begin{definition}[Measured Software Scale-up]
Let 
\begin{equation}\label{measurscale}
Sc_{1,QP}^{meas}(\mathcal{A^{DD}}):= \frac{T_{flop}(N_{loc})}{QP \cdot (T_{flop}(N_{loc})+T_{oh}(N_{loc}))}
\end{equation}
be the measured software scale-up in going from $1$ to $QP$.
\begin{flushright}
$\spadesuit$
\end{flushright}

\end{definition}

\begin{proposition} Let $s_{nproc}^{loc}(\mathcal{A}^{loc}_{4DVar})$ denote the speedup of the local  parallel algorithm $(\mathcal{A}^{loc}_{4DVar})$. If $$0 \leq \frac{S}{V}(\mathcal{A}^{loc}_{4DVar})<1-\frac{1}{s_{nproc}^{loc}(\mathcal{A}^{loc}_{4DVar})}\quad ,$$ then it holds that
\begin{equation}\label{Eq_relScale}
 Sc_{1,QP}^{meas}(\mathcal{A}^{DD}_{4DVar})= \alpha(N_{loc},QP) Sc_{1,QP}(\mathcal{A}^{DD}_{4DVar})
\end{equation}
with
\begin{eqnarray}\label{alpha}
\alpha(N_{loc},QP)(\mathcal{A}^{DD}_{4DVar}) & = &  \frac{T_{flop}(N_{loc})}{\frac{QP \,T_{flop}(N_{loc})}{s_{nproc}^{loc}(\mathcal{A}^{loc}_{4DVar})}+QP\,T_{oh}(N_{loc})}\nonumber \\
& = &\frac{s_{nproc}^{loc}(\mathcal{A}^{loc}_{4DVar})\frac{T_{flop}(N_{loc})}{QP\,T_{flop}(N_{loc})}}{1+\frac{s_{nproc}^{loc}(\mathcal{A}^{loc}_{4DVar})T_{oh}(N_{loc})}{T_{flop}(N_{loc})}}.
\end{eqnarray}

\noindent If  $$\alpha(N_{loc},QP):= \frac{s_{nproc}^{loc}(\mathcal{A}^{loc}_{4DVar})}{1+\frac{s_{nproc}^{loc}(\mathcal{A}^{loc}_{4DVar})T_{oh}(N_{loc})}{T_{flop}(N_{loc})}}= \frac{s_{nproc}^{loc}(\mathcal{A}^{loc}_{4DVar})}{1+s_{nproc}^{loc}(\mathcal{A}^{loc}_{4DVar})\frac{S}{V}(\mathcal{A}^{loc}_{4DVar})}$$
from (\ref{alpha}), it becomes the thesis in (\ref{Eq_relScale}).

\begin{flushright}
$\clubsuit$
\end{flushright}

\end{proposition}
\noindent In the following we  denote the measured scale-up as $Sc_{1,QP}^{meas}(\mathcal{A}^{DD}_{4DVar})$ or as $Sc_{1,QP}^{meas}(N)$, respectively.\\

\noindent The next proposition allows us to examine the benefit on the measured  scale-up arising  from the speedup of the local  parallel algorithm $s_{nproc}^{loc}(\mathcal{A}^{loc}_{4DVar})$, mainly in the presence of a multilevel decomposition, where $s_{nproc}^{loc}(\mathcal{A}^{loc}_{4DVar}) >1$.

\begin{proposition} It holds that
$$s_{nproc}^{loc}(\mathcal{A}^{loc}_{4DVar})  \in [1,QP]\Rightarrow Sc_{QP}^{meas}(\mathcal{A}^{DD}_{4DVar})  \in ]Sc_{1,QP}(\mathcal{A}^{DD}_{4DVar}) , QP \,Sc_{1,QP}(\mathcal{A}^{DD}_{4DVar}) [.$$
\noindent {\bf Proof:}\\
\begin{itemize}
  \item  If $s_{nproc}^{loc}(\mathcal{A}^{loc}_{4DVar}) =1$, then $$\alpha(N,QP) <1 \Leftrightarrow Sc_{1,QP}^{meas}(\mathcal{A}^{DD}_{4DVar}) < Sc_{1,QP}(\mathcal{A}^{DD}_{4DVar}) .$$
  \item If $s_{nproc}^{loc}(\mathcal{A}^{loc}_{4DVar}) >1$,  then
$$  \alpha(N,QP) >1   \Leftrightarrow  Sc_{1,QP}^{meas}(\mathcal{A}^{DD}_{4DVar})  > Sc^f_{1,QP}(\mathcal{A}^{DD}_{4DVar}) .$$
\item If $s_{nproc}^{loc} (\mathcal{A}^{loc}_{4DVar})=QP$, then $$  1< \alpha(N,QP) <QP   \Rightarrow   Sc_{1,QP}^{meas}(\mathcal{A}^{DD}_{4DVar})  < QP \cdot Sc^f_{1,QP}(\mathcal{A}^{DD}_{4DVar}) .$$
\end{itemize}

\begin{flushright}
$\clubsuit$
\end{flushright}

\end{proposition}

\noindent We  may conclude the following:
 \begin{enumerate}
   \item Strong scaling: if $QP$ increases and $M\times K$ is fixed, the scale-up factor increases but  the surface-to-volume ratio also increases.
   \item Weak scaling: if $QP$ is fixed and $M \times K$ increases, the scale-up factor stagnates and the surface-to-volume ratio decreases.
 \end{enumerate}

 \noindent  Thus, one needs to find the appropriate value of the number of subdomains, $QP$, giving the right tradeoff between the scale-up and the overhead of the algorithm.

\section{Scalability results}
\noindent 
 The results presented here are just a starting point toward  the assessment of the software scalability. More precisely, we introduce   simplifications and assumptions appropriate for a proof-of-concept study  in order to  get   values of the measured scale-up  of the one iteration of the  parallel algorithm.\\ 
 \noindent Since the main outcome of the decomposition is that the parallel algorithm is oriented to better exploit the high performance of new architectures where concurrency is implemented both at the coarsest and finest levels of granularity, such as a distributed-memory multiprocessor (MIMD) and a graphics processing unit (GPU),  we consider a distributed-computing environment located in the University of Naples Federico II campus, connected by local-area network made of the following:\\

\begin{itemize}
 \item $PE_1$ (for the coarsest level of granularity): a MIMD architecture made of $8$ nodes that consist of distributed-memory DELL M600 blades connected by a 10 Gigabit Ethernet technology. Each blade consists of $2$ Intel Xeon@2.33GHz quadcore processors sharing the same local 16 GB of RAM memory for a total of $8$ cores per blade and  $64$ total cores.   \\
\item $PE_2$ (for the finest level of granularity): a Kepler architecture of the GK110 GPU  \cite{k20},  which consists of a set of 13 programmable single-instruction, multiple-data (SIMD) streaming multiprocessors (SMXs), connected to a quad-core Intel i7 CPU running at 3.07 GHz, 12 GB of RAM. For host(CPU)-to-device(GPU) memory transfers CUDA-enabled graphic cards are connected to a PC motherboard via a PCI-Express (PCIe) bus \cite{pcie}.  
For this architecture the maximum number of active threads per multiprocessor is 2,048, which means that the maximum number of active warps per SMX is 64.
\end{itemize}

\noindent  Our implementation uses the matrix and vector functions in the Basic Linear Algebra Subroutines (BLAS) for $PE_1$ and the CUDA Basic Linear Algebra Subroutines (CUBLAS) library for $PE_2$. The routines used for computing the minimum of $J$ on $PE_1$ and $PE_2$ are described in \cite{lbfgsNoc} and \cite{lbfgsgpu}, respectively.\\
\noindent The case study is based on the shallow water equations  on the sphere.
The SWEs have been used extensively as a simple model of the atmosphere or ocean circulation because they
contain the essential wave propagation mechanisms found in general circulation models \cite{ShallowWater}.\\

\noindent The SWEs in spherical coordinates are

\begin{eqnarray}
\frac{\partial u}{\partial t} & = & - \frac{1}{a \cos{\theta}} \left( u \frac{\partial u}{\partial \lambda} + v \cos{\theta}\frac{\partial u}{\partial\theta}\right)
                                    + \left( f + \frac{u\tan{\theta}}{a}\right) v
                                    - \frac{g}{a \cos{\theta}} \frac{\partial h}{\partial \lambda} \label{EquU}\\
\frac{\partial v}{\partial t} & = & - \frac{1}{a \cos{\theta}} \left( u \frac{\partial v}{\partial \lambda} + v \cos{\theta}\frac{\partial v}{\partial\theta}\right)
                                    + \left( f + \frac{u\tan{\theta}}{a}\right) u
                                    - \frac{g}{a} \frac{\partial h}{\partial \theta} \label{EquV}\\
\frac{\partial h}{\partial t} & = & - \frac{1}{a \cos{\theta}} \left(\frac{\partial\left(hu\right)}{\partial \lambda} + \frac{\partial\left(hu\cos{\theta}\right)}{\partial \theta}.\right) \label{EquH}
\end{eqnarray}

\noindent Here $f$ is the Coriolis parameter given by $f = 2\Omega\sin{\theta}$, where $\Omega$ is the angular speed of the rotation of the Earth;
$h$ is the height of the homogeneous atmosphere (or of the free ocean surface); $u$ and $v$ are the zonal and meridional wind (or the ocean velocity)
components, respectively; $\theta$ and $\lambda$ are the latitudinal and longitudinal directions, respectively; and $a$ is the
radius of the Earth and $g$ is the gravitational constant.

\noindent We  express the system of equations (\ref{EquU})--(\ref{EquH}) using a compact form:

\begin{equation}
\frac{\partial{\mathbf{Z} }}{\partial t} = \mathcal{M}_{t-\Delta t\rightarrow t}\left({\mathbf{Z}}, \right) \label{ShallowWaterEqu}
\end{equation}
where
\begin{equation}
{\mathbf{Z}}=\left(
\begin{array}{c}        
       u \\ v \\ h
\end{array}     
  \right)
\end{equation}
and

\begin{eqnarray}
\mathcal{M}_{t-\Delta t\rightarrow t}\left({\mathbf{Z}} \right) & = & \left(
                        \begin{array}{l}
- \frac{1}{a \cos{\theta}} \left( u \frac{\partial u}{\partial \lambda} + v \cos{\theta}\frac{\partial u}{\partial\theta}\right)
                                    + \left( f + \frac{u\tan{\theta}}{a}\right) v
                                    - \frac{g}{a \cos{\theta}} \frac{\partial h}{\partial \lambda} \\
- \frac{1}{a \cos{\theta}} \left( u \frac{\partial v}{\partial \lambda} + v \cos{\theta}\frac{\partial v}{\partial\theta}\right)
                                    + \left( f + \frac{u\tan{\theta}}{a}\right) u
                                    - \frac{g}{a} \frac{\partial h}{\partial \theta} \\
- \frac{1}{a \cos{\theta}} \left(\frac{\partial\left(hu\right)}{\partial \lambda} + \frac{\partial\left(hu\cos{\theta}\right)}{\partial \theta}\right)
                        \end{array}
                         \right) \nonumber \\
                      & = & \left(
                       \begin{array}{l}
                             F_1 \\
                             F_2 \\
                             F_3
                       \end{array}
                         \right).
\end{eqnarray}

\noindent We discretize (\ref{ShallowWaterEqu}) just in space using an unstaggered Turkel--Zwas scheme \cite{Turkel-Zwas1,Turkel-Zwas2}, and
we obtain
\begin{equation}
\frac{\partial{\mathbf{Z}_{disc}}}{\partial t} = \mathcal{M}^{t-\Delta t\rightarrow t}_{disc}\left({\mathbf{Z}_{disc}}, \right) \label{ShallowWaterSemiDiscrete}
\end{equation}
where
\begin{equation}
{ \mathbf{Z}_{disc}}=\left(
\begin{array}{c}        
       \left(u_{i,j}\right)_{i=0,\ldots,nlon-1;j=0,\ldots,nlat-1} \\
       \left(v_{i,j}\right)_{i=0,\ldots,nlon-1;j=0,\ldots,nlat-1} \\
       \left(h_{i,j}\right)_{i=0,\ldots,nlon-1;j=0,\ldots,nlat-1} \\
\end{array}     
  \right)
\end{equation}
and
\begin{equation}
\mathcal{M}^{t-\Delta t\rightarrow t}_{disc}\left({\mathbf{Z}_{disc}} \right) = \left(
                        \begin{array}{l}
       \left(U_{i,j}\right)_{i=0,\ldots,nlon-1;j=0,\ldots,nlat-1} \\
       \left(V_{i,j}\right)_{i=0,\ldots,nlon-1;j=0,\ldots,nlat-1} \\
       \left(H_{i,j}\right)_{i=0,\ldots,nlon-1;j=0,\ldots,nlat-1}. \\
                        \end{array}
                         \right)
\end{equation}
Thus 
\begin{eqnarray*}
U_{i,j}  & = & - \sigma_{lon} \frac{u_{i,j}}{\cos{\theta_j}} \left( u_{i+1,j} - u_{i-1,j}\right)  \\
         &   & - \sigma_{lat} \ {v_{i,j}}                \left( u_{i,j+1} - u_{i,j-1}\right) \\
         &   & - \sigma_{lon} \frac{g}{p \cos{\theta_j}} \left( h_{i+p,j} - h_{i-p,j}\right) \\
         &   & + 2 \left[ \left(1-\alpha\right) \left( 2\Omega\sin{\theta_j}+ \frac{u_{i,j}}{a}\tan\theta_j\right)v_{i,j}     \right. \\
         &   & +   \left. \frac{\alpha}{2}      \left( 2\Omega\sin{\theta_j}+ \frac{u_{i+p,j}}{a}\tan\theta_j\right)v_{i+p,j} \right. \\
         &   & +   \left. \frac{\alpha}{2}      \left( 2\Omega\sin{\theta_j}+ \frac{u_{i-p,j}}{a}\tan\theta_j\right)v_{i-p,j} \right] \\
V_{i,j}  & = & - \sigma_{lon} \frac{u_{i,j}}{\cos{\theta_j}} \left( v_{i+1,j} - v_{i-1,j}\right)  \\
         &   & - \sigma_{lat} \ {v_{i,j}} \left( u_{i,j+1} - u_{i,j-1}\right) \\
         &   & - \sigma_{lat} \frac{g}{q} \left( h_{i,j+q} - h_{i,j-q}\right) \\
         &   & - 2 \left[ \left(1-\alpha\right) \left( 2\Omega\sin{\theta_j}+ \frac{u_{i,j}}{a}\tan\theta_j\right)u_{i,j}             \right. \\
         &   & +   \left. \frac{\alpha}{2}      \left( 2\Omega\sin{\theta_{j+q}}+ \frac{u_{i,j+q}}{a}\tan\theta_{j+q}\right)u_{i,j+q} \right. \\
         &   & +   \left. \frac{\alpha}{2}      \left( 2\Omega\sin{\theta_{j-q}}+ \frac{u_{i,j-q}}{a}\tan\theta_{j-q}\right)u_{i,j-q} \right] \\
H_{i,j}  & = & - \alpha \left\{ \frac{u_{i,j}}{\cos{\theta_j}}  \left( h_{i+1,j} - h_{i-1,j}\right) \right.\\
         &   & + \ {v_{i,j}} \left( h_{i,j+1} - h_{i,j-1}\right) \\
         &   & + \frac{h_{i,j}}{\cos{\theta_j}} \left[ \left(1-\alpha\right) \left( u_{i+p,j} - u_{i-p,j}\right)  \right. \\
         &   & + \left. \frac{\alpha}{2} \left(u_{i+p,j+q} - u_{i-p,j+q} + u_{i+p,j-q} - u_{i-p,j-q}\right)\right]\frac{1}{p} \\
         &   & + \left[ \left(1-\alpha\right) \left( v_{i,j+q}\cos{\theta_{j+q}} - v_{i,j-q}\cos{\theta_{j-q}}\right) \right. \\
         &   & + \frac{\alpha}{2} \left( v_{i+p,j+q}\cos{\theta_{j+q}} - v_{i+p,j-q}\cos{\theta_{j-q}} \right)\\
         &   &\left.\left.+\frac{\alpha}{2} \left( v_{i-p,j+q}\cos{\theta_{j+q}} - v_{i-p,j-q}\cos{\theta_{j-q}} \right)\right]\frac{1}{q}. \right\}
\end{eqnarray*}



\noindent The numerical model depends on a combination  physical parameters, including the number of state variables in the
model, the number of observations in an assimilation
cycle, and the numerical parameters as the discretization step in time and  in space  are defined on the basis of a discretization grid used by data available in the  Ocean Synthesis/Reanalysis Directory of Hamburg University (\cite{Dati}). 

\noindent Our data assimilation experiments are initialized by choosing snapshots from the run prior to the start of the assimilation experiment and treating it as realization valid at the nominal  time. Then, the model state is advanced to the next time using the forecast model, and the observations are combined with the forecasts (i.e., the background) to produce the analysis. This process is iterated. As it proceeds, the process fills gaps in sparsely observed regions, converts observations to improved estimates of model variables, and filters observation noise. All this is done in a manner that is physically consistent with the dynamics of the ocean as represented by the model. In our experiments, the simulated observations are created by sampling the model states  and adding random errors to those values. A detailed description of the simulation, together with  the results and the software implemented,  is presented in \cite{arxiv}. In the following, we focus mainly on performance results.  \\

\noindent The reference domain decomposition strategy uses the following correspondence between $QP$  and $nproc$,
$$QP\leftrightarrow nproc,$$

\noindent which means that the number of subdomains coincides with the number of available processors.  \\
\noindent According to the characteristics of the physical domain in SWEs,  the total number of grid points in space   is  $$M=nlon \times nlat\times n_z \quad .$$ Assume that   $$nlon=nlat=n,$$  where $n_z=3$.  Since the unknown vectors are the fluid height or depth and the two-dimensional fluid velocity fields,  the problem size in space  is $$M=n^2 \times 3\,.$$ We assume a 2D uniform domain decomposition along the latitude-longitude directions such that  

\begin{equation}\label{Eq-r}
D_s:= \frac{M}{p}=  nloc_x \times  nloc_y \times 3
\end{equation}

\noindent with

\begin{equation}\label{nlocx}
nloc_x := \frac{n }{p_1} +2o_x\,\,,\, nloc_y := \frac{n}{p_2} +2o_y\,\,\,,\,n_z := 3\,\,,
\end{equation}
where $p_1 \times p_2= p$. Here $o_x$ and $o_y$ denote the overlapping regions along $x$ and $y$ directions.

\noindent  Since the GPU ($PE_2$) can process only the data in its global memory, in a generic parallel algorithm execution the host acquires this input data and sends it to the device memory, which concurrently calculates the minimization of the 4D-Var functional. To avoid continuous relatively slow data transfer from the host to the device and  to reduce the overhead, we  store the device with the entire work data prior to any processing. Specifically,  the maximum value of $D_s$ in (\ref{Eq-r}) is chosen such that the amount of data related each subdomain (we denote it with $Data_{mem} (Mbyte)$) can be completely stored in the memory.\\

\noindent If we assume that $nloc_x=nloc_y$ and we let $n_{loc}=nloc_x=nloc_y$, since the global GPU memory is  5 GB, we have the values of usable $n_{loc}$ described in Table \ref{Tab:space}, Table \ref{fig:fig5} reports the values of the speedup $s_{nproc}^{loc}$ in terms of gain obtained by using the GPU versus the CPU. We note that  CUBLAS routines allow us to reduce on average  18 times the execution time necessary for a single CPU for the minimization part.\\

\begin{table}[!ht]
\centering
\begin{tabular}{c|cccccccccc}
\hline
    $n_{loc}$                           & $32 $ & $40 $  & $48 $ & $56 $   & $64 $ & $72 $   & $80 $   & $88 $ \\
    
															\hline\\

 $Data_{mem} (Mbyte)$        & 177   & 286   & 485  & 812    & 1313 & 2041   & 3057   & 4427 \\
 
\hline
\end{tabular}
\caption{The amount of memory required to store data related to each subdomain  on $PE_2$ expressed in Mbyte.}\label{Tab:space}
\end{table}


\begin{table}[!ht]
\centering
\begin{tabular}{|c|cccccccc}
\hline
  $  n_{loc}   $                       &   32  & 40 & 48 & 56&   64 & 72   & 80   & 88 \\ \hline
$\frac{T_{blas}}{T_{cublas}}$         & 15.3& 17.5 & 18.08 & 19.0 & 19.8 & 20.2 & 22.5& 20.54 \\[.1cm]
\hline
\end{tabular}
\caption{Values of the speedup $s_{nproc}^{loc}$ in terms of gain obtained by using the GPU versus the CPU. The CUBLAS routines allow  reducing on average by 18 times the execution time necessary for a single CPU for the minimization part.
}
\label{fig:fig5}
\end{table}

         





\begin{figure}
    \centering
    \includegraphics[width=8cm]{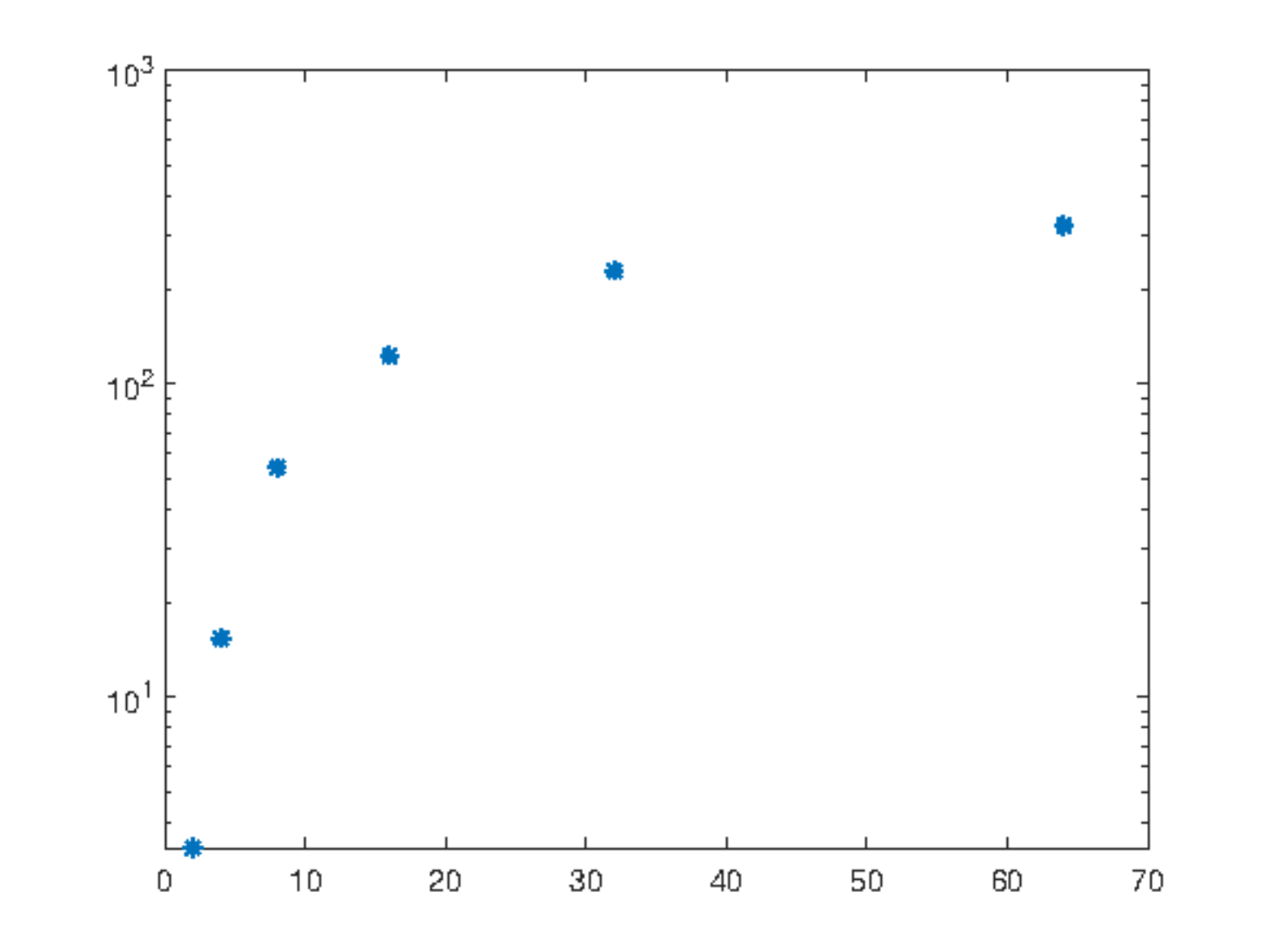}
    \caption{{Weak scalability   of one iteration of the  parallel algorithm $\mathcal{A}^{DD}_{4DVar}$  with   $n_{loc}=32$ computed by sing the measured software scale-up $Sc^{meas}_{1,QP}$ defined in (\ref{measurscale}).  }\label{my-label2f}}
    \label{fig:my_label}
\end{figure}

\begin{table}[!ht]
\centering
\begin{tabular}{c|cccccc}
\hline
         $QP$  & 2     & 4 & 8   & 16 &32 &64  \\

         \hline\\											

        problem size       & $6.1\cdot 10^3$ & $1.2\times 10^4$ & $2.4\cdot 10^4$ &$4.9\cdot 10^4$ &$9.8\cdot 10^4$&$1.9 \times 10^5$  \\

         \hline\\

$Sc^{meas}_{1,QP}$   & $3.3\cdot 10^0$ & $1.54\cdot 10^1$   & $5.41\cdot 10^1$ &  $1.23\cdot 10^2$ & $2.30\cdot 10^2$ & $3.2 \times 10^2$\\
\hline
\end{tabular}
\caption{Weak scalability   of one iteration of the  parallel algorithm $\mathcal{A}^{DD}_{4DVar}$  with   $n_{loc}=32$ computed by using the measured software scale-up $Sc^{meas}_{1,QP}$ defined in (\ref{measurscale}).}\label{my-label2}
\end{table}
\noindent The outcome  from these experiments is that the algorithm scales up according to the performance analysis (see Figure \ref{my-label2f}). Indeed, as expected, as $QP$ increases, the scale-up factor increases and  the surface-to-volume ratio increases, too, so that performance gain tends to become  stationary. This the inherent tradeoff between speedup and efficiency of any software architecture.   


\section{Conclusions}
 We provide a  complete computational   framework of a space-time decomposition approach for 4D-Var. This includes the mathematical framework, the numerical algorithm, and  its performance validation. We measure the performance of the  algorithm using a simulation case study based on the SWEs on the sphere. Results presented here are just a starting point toward  the assessment of the software scalability. More precisely, we introduce   simplifications and assumptions appropriate for a proof-of-concept study  in order to  measure scale-up  of  one iteration of the  parallel algorithm.  The overall insight we get from these experiments is that the algorithm scales up according to the performance analysis. 
 
\noindent We are currently working on the development of a flexible framework    ensuring efficiency and code readability,  exploiting future technologies, and  including a quantitative assessment of scalability. In this regard, we could combine the proposed approach with the PFASST algorithm. Indeed, PFASST could be concurrently employed as a local solver of each reduced-space PDE-constrained optimization subproblem, exposing even more temporal parallelism.
\noindent This framework will allow  designing,  planning, and running  simulations to identify and overcome the limits of this approach.   
\section*{Acknowledgments} This work was developed within the research activity of the H2020-MSCA-RISE-2016
554 NASDAC Project N. 691184. This work has been realized thanks to the use of the S.Co.P.E. computing infrastructure at the University of Naples.
The material is based upon work supported by the U.S. Department of Energy, Office of Science, under contract DE-AC02-06CH11357.

\section{Declarations}
The authors confirm that the research described in this  work has not received any funds.\\
\noindent The authors confirm that there are not any conflicts of interest.\\
\noindent The authors confirm that data and code can be available at request.


\end{document}